\documentclass[a4paper]{article}

\usepackage[english]{babel}
\usepackage[utf8x]{inputenc}
\usepackage{amsmath}
\usepackage{graphicx}
\usepackage[colorinlistoftodos]{todonotes}
\usepackage{authblk}

\usepackage{amsmath}
\usepackage{epsfig,amsthm}
\usepackage{latexsym}
\usepackage{amsfonts}
\usepackage{amssymb}
\usepackage{amscd}
\usepackage{mathrsfs}
\usepackage[all,cmtip]{xy}
\usepackage{enumerate}
\usepackage{tikz}
\usepackage{arydshln}
\usepackage{bbm}
\usepackage{mathtools}
\usepackage{leftindex}
\usepackage{mathdots}
\usepackage{centernot}

\usepackage{tocloft}
\usepackage{amsfonts}
\usepackage{color}
\usepackage[colorlinks]{hyperref}
\hypersetup{pdfstartview={FitH},         linkcolor=blue,  citecolor=green }

\numberwithin{equation}{section} 

\usepackage[a4paper,bindingoffset=0.2in,left=1in,right=1in,top=1in,bottom=1in,footskip=.25in]{geometry}

\usepackage{fancyhdr}
\AtEndDocument{%
  \par
  \medskip
  \begin{tabular}{@{}l@{}}%
  \textsc{Radboud University}
  \\
    \textsc{Heyendaalseweg 135
, 6525AJ Nijmegen, the Netherlands
}\\
    \textit{E-mail address}: \texttt{\href{mailto:geotam@science.ru.nl}{geotam@science.ru.nl}}
  \end{tabular}}

\newtheorem{thm}{Theorem}[section]
\newtheorem{cor}[thm]{Corollary}
\newtheorem{lem}[thm]{Lemma}
\newtheorem{prop}[thm]{Proposition}

\theoremstyle{definition}

\newtheorem{rmk}[thm]{Remark}

\newcommand{\Fo}{{F_{\bullet}}}
\newcommand{\Eo}{{E_{\bullet}}}
\newcommand{\qo}{{q_{\bullet}}}

\newcommand{\tr}{{\mathrm{tr}}}

\newcommand{\GL}{{\mathrm{GL}}}

\newcommand{\SL}{{\mathrm{SL}}}
\newcommand{\SO}{{\mathrm{SO}}}
\newcommand{\SP}{{\mathrm{Sp}}}

\newcommand{\sgn}{{\mathrm{sgn}}}

\newcommand{\Ind}{{\mathrm{Ind}}}
\newcommand{\cInd}{{\mathrm{cInd}}}
\newcommand{\Res}{{\mathrm{Res}}}

\newcommand{\St}{{\mathrm{St}}}

\newcommand{\diag}{{\mathrm{diag}}}
\newcommand{\antidiag}{{\mathrm{anti}\text{-}\mathrm{diag}}}

\newcommand{\Hom}{{\mathrm{Hom}}}

\newcommand{\Gal}{{\mathrm{Gal}}}

\newcommand{\Ad}{{\mathrm{Ad}}}

\newcommand{\Red}{{\mathrm{Red}}}

\newcommand*{\thead}[1]{\multicolumn{1}{|c|}{\bfseries #1}}

\title{Depth Zero Supercuspidal Representations of Classical Groups into L-packets: the Typically Almost Symmetric Case}
\author{
 Geo Kam-Fai Tam }
 \affil{ 
 \href{mailto:geotam@science.ru.nl}{geotam@science.ru.nl}
}
\date{\vspace{-0.5cm}\today}

\begin{document}
\maketitle

\raggedbottom
\interlinepenalty=10000

\begin{abstract}
We classify what we call ``typically almost symmetric'' depth zero supercuspidal representations of classical groups into L-packets. Our main results resolve an ambiguity in the paper of Lust-Stevens \cite{Lust-Stevens} in this case, where they could only classify these representations in two or four, if not one, L-packets. By assuming the expected numbers of supercuspidal representations in the L-packets, we employ only simple properties of the representations to prove the main results. In particular, we do not require any deep calculations of character values.
\end{abstract}

\tableofcontents

\section{Introduction}

This paper is an attempt to resolve the ambiguities in \cite{Lust-Stevens} about classifying depth zero supercuspidal representations of classical groups over a local field into unique L-packets, under the local Langlands correspondence (LLC) for these groups. Our attempt is successful in what we call the \emph{typically almost symmetric} case, which will be explained shortly.

Recall that, in the classification of the mentioned paper, we aim at describing the Langlands parameter, or equivalently, the Jordan blocks of depth zero supercuspidal representations in terms of their inducing cuspidal types. The description should agree with the theory of endoscopy asserted by 
\cite{Arthur-new-book, Arthur-even-orthogonal, Mok-unitary, unitary-inner}, an important branch of the Langlands program. This description has been successful for unramified unitary groups \cite[Sec 9.2]{Lust-Stevens}, and in some cases for other classical groups, but in general we can only classify such a representation in a union of one, two or four packets, due to an ambiguity from the twists by unramified characters of the components of the prescribed parameters. Also, it seems to the author that the descriptive results in \cite{Lust-Stevens} are only proved in principle: a precise description of the correspondence seems to be not bookkept in \emph{loc. cit.} In particular, we want to know whether the Langlands parameter of a depth zero self-dual supercuspidal representation of a general linear group has orthogonal or symplectic image. 

As mentioned in the abstract, we resolve these issues for depth zero typically almost symmetric supercuspidal representations. The precise results are given in Section \ref{subsection The main result}, and summarized in Theorem \ref{intro main theorem} below. Before stating the theorem, let's emphasize one {\bf special feature} about the proof of our results: beyond the foundational knowledge on covering types and Hecke algebras for depth zero supercuspidal representations, which are already covered in \cite{Lust-Stevens} and other references, we only require a few additional simple properties of these representations; the knowledge of their exact character values are not needed. In general, computing these values requires advanced techniques such as character sheaves (see \cite{Waldspurger-Lusztig-conj} for example), but we do not need them in our case.

To state the theorem, we first require some knowledge on the inducing cuspidal types for constructing the mentioned representations. Recall from \cite{Moy-Prasad-Jacquet} that a depth zero supercuspidal representation $\pi$ of a connected reductive group $G$ over a non-Archimedean field $F$ is compactly induced from an extension, to the normalizer $\mathcal J$ of a maximal parahoric subgroup $\mathcal J_0$, of the inflation of a cuspidal representation $\rho$ of $\mathcal J_0$, colloquially called a type for $\pi$. When $G$ is a classical group (defined as the fixed-point subgroup of an involution $\sigma$ on a general linear group), the reductive quotient of $\mathcal J_0$ is the connected component of a product of at most two classical groups over the residue field $\mathbb F_q$ of $F$, denoted by $\mathsf {G}_y$ and $\mathsf {G}_z$ in this paper, which are possibly disconnected when they are of orthogonal type. In fact, this observation encourages us to slightly expand our scope to also consider $G$ to be a (disconnected) full even orthogonal group, so $\mathcal J$ is indeed a product of at most two orthogonal groups over $\mathbb F_q$.

We now assume that $G$ is even orthogonal, as well as each $\mathsf {G}_w$, for $w\in \{y,z\}$. The classification of cuspidal representations of $\mathsf {G}_w$ is known due to Lusztig \cite{lusz-classical}. Here we use a modification by Waldspurger \cite{Waldspurger-Lusztig-conj}, and denote by 
$$\rho(\epsilon_+  m_+,\epsilon_- m_-,(m_{\tilde P})_{\tilde P\neq X\pm 1}), \quad \epsilon_+,\epsilon_-\in \{\pm 1\},$$
the cuspidal representations of $\mathsf {G}_w$ lying in the Lusztig series $\mathcal E(\mathsf {G}_w,s)$, where the numbers $m_+=m_{X-1}$, $m_-=m_{X+1}$, and other $m_{\tilde P}$ all belong to $\mathbb N_{\geq 0}$, and $s$ is a semisimple element in the dual group of  $\mathsf {G}_w$ whose characteristic polynomial (by viewing $\mathsf {G}_w$ as a matrix group) is 
$$(X-1)^{2m_+}(X+1)^{2m_-}\prod_{
\tilde{P}}\tilde P^{m_{\tilde P}(m_{\tilde P}+1)/2},$$
with  $\tilde{P} $ ranging over all irreducible polynomials in $\mathbb{F}_q[X]$ of degree $>1$ and which are self-dual (for any root in $\tilde{P}$ its inverse is also a root).

We can now define our condition on depth zero supercuspidal representations: we call $\pi$ \emph{typically almost symmetric} if its inducing type $\rho = \rho_y\times \rho_z$ on $\mathsf G_y\times \mathsf G_z$ satisfies $|m_{\tilde P,y}| = |m_{\tilde P,z}|$ for all self-dual irreducible polynomials $\tilde P$. Therefore, with fixed data $( |m_+|,|  m_-|,(m_{\tilde P})_{\tilde P\neq X\pm 1})$, we can produce at most four inducing types: they are called \emph{companions} of each other.

Note that each $\tilde P$ corresponds to a cuspidal representation of $\GL(m,\mathbb F_q)$ where $m=\deg \tilde P$. If $\tilde{\rho}$ is such a representation, the inflation of $\tilde\rho$ to $\tilde{\mathcal J}_0:=\GL(m,\mathcal O_F)$ admits two self-dual extensions, denoted by $\tilde\rho$ and $\tilde\rho'$, to the normalizer $\tilde {\mathcal J} := F^\times \GL(m,\mathcal O_F)$ of $\tilde{\mathcal J}_0$ in $\tilde{G} = \GL(m,F)$. Then $\tilde\rho$ and $\tilde\rho'$ respectively compactly induce to two non-isomorphic depth zero supercuspidal representations $\tilde\pi$ and $\tilde\pi'$ of  $\tilde{G}$.

On the Galois side, a Langlands parameter $\varphi$ for $G$ is viewed as an isomorphism class of a self-dual  representation of the Weil-Deligne group $\mathcal W_F\times \SL(2,\mathbb C)$ with an even degree (determined by $G$) and an orthogonal image. To describe a depth zero parameter admitting cuspidal representations in its packet, we first denote by $\St[a]$ the degree-$a$ irreducible algebraic representation of $\SL(2,\mathbb C)$, and put
$\St[[a]] = \St[a] \oplus \St[a-2]\oplus \cdots$, where $\St[a]$ is the zero representation for all $a\in -\mathbb{N}_{\geq 0}$. Then the L-packet of $\varphi$ contains a supercuspidal representation if it takes the form
$$\varphi = 1\otimes \St[[a]] \oplus 
\omega_0\otimes \St[[b]] \oplus 
\omega_1\otimes \St[[c]] \oplus 
\omega_2\otimes \St[[d]] \oplus 
\bigoplus_{ 
\tilde\varphi }\tilde\varphi\otimes \St[[a_{\tilde\varphi}]].$$
Here $1,\omega_0,\omega_1,\omega_2$ are the four quadratic characters of $\mathcal W_F$: the trivial, the unramified quadratic, and the two ramified quadratic with $\omega_j(\varpi) = (-1)^{j-1}$ for $j\in \{1,2\}$, the numbers $a,b,c,d$ are odd with certain conditions on their mod-4 values, and $\tilde\varphi $ ranges over self-dual depth zero irreducible representations of  $\mathcal W_F$ of degree $>1$, whose image lies in an orthogonal (resp. symplectic) group if and only if the attached integer $a_{\tilde\varphi}$ is odd (resp. even).

Note that if the image of $\tilde\varphi$ is orthogonal, then its unique self-dual non-equivalent twist $\tilde\varphi'$ by a suitable unramified character is symplectic. Suppose that $\tilde\varphi$ and $\tilde\varphi'$ correspond to $\tilde\pi = \tilde\pi_{\tilde\varphi}$ and $\tilde\pi' = \tilde\pi_{\tilde\varphi'}$ under the LLC for general linear groups, then one can establish a bijection
$\{\tilde\varphi,\tilde\varphi'\} \leftrightarrow \{\tilde P\}$ between such a pair of parameters and self-dual irreducible polynomials over $\mathbb F_q$.

For symplectic and unramified unitary groups, we refer to Sections \ref{subsection Finite classical groups} and \ref{section LLC for classical groups} for similar definitions and constructions of their representations and parameters, together with modifications. We now provide our main results, summarizing Propositions \ref{main-result-symplectic}, \ref{main-result-orthogonal}, \ref{main result for unramified unitary}, and Corollary \ref{Corollary for parity on representation side}.

\begin{thm} 
\label{intro main theorem}
\begin{enumerate}[(i)]
\item The following tables list the parameter $\varphi$ whose L-packet contains the  supercuspidal representations $\pi = \cInd_{\mathcal J}^G \rho$ of the listed classical group $G$, where the inducing type $\rho= \rho_y\times \rho_z$ are inflated from the listed cuspidal representations.

$G=\SP(4n,F)$, $\mathsf G_y = \mathsf G_z = \SP(2n,\mathbb F_q)$ \vspace{-5pt}
\renewcommand{\arraystretch}{1.4}
\begin{center}
\begin{tabular}{|l |c |} 
 \hline
  \thead{ $\varphi$} & $\{\rho_y,\rho_z\}$ 
\\ 
 \hline
 $1\otimes \St[[4m_++1]]\oplus \omega_1\otimes \St[[4m_--1]]$ 
& 
 $\rho_y = \rho_z = \rho(m_+, \pm m_-,(m_{\tilde P})_{\tilde P})$
 \\
\qquad\qquad\qquad\qquad\qquad\qquad $\bigoplus_{\tilde\varphi }\tilde\varphi \otimes \St[[2m_{\tilde\varphi}]]$ &
\\ 
 $1\otimes \St[[4m_++1]]\oplus \omega_2\otimes \St[[4m_--1]]$
& 
$\{\rho( m_+,  m_-,(m_{\tilde P})_{\tilde P}),\rho( m_+, - m_-,(m_{\tilde P})_{\tilde P})\}$
\\
\qquad\qquad\qquad\qquad\qquad\qquad 
$\bigoplus_{\tilde\varphi }\tilde\varphi \otimes \St[[2m_{\tilde\varphi}]]$ 
&
\\
 \hline
\end{tabular}
\end{center}

$G=\mathrm O(4n,F)$, $\mathsf G_y = \mathsf G_z = \mathrm {O}(2n,\mathbb F_q)$ \vspace{-15pt}
\begin{center}
\begin{tabular}{|l |c |} 
 \hline
\thead{ $\varphi$ } & $\{\rho_y,\rho_z\}\text{ with }{\epsilon_+, \epsilon_-}\in \{\pm 1\} $ 
\\ 
 \hline
 $1\otimes \St[[4m_+-1]]\oplus \omega_1\otimes \St[[4m_--1]]$ 
& 
 $\rho_y = \rho_z = \rho(\epsilon_+ m_+, \epsilon_- m_-,(m_{\tilde P})_{\tilde P})$
  \\
\qquad\qquad\qquad\qquad\qquad\qquad $\bigoplus_{\tilde\varphi }\tilde\varphi \otimes \St[[2m_{\tilde\varphi}]]$ &
\\ 
 $\omega_0\otimes \St[[4m_+-1]]\oplus \omega_1\otimes \St[[4m_--1]]
 $ 
& 
$\{\rho(\epsilon_+ m_+, \epsilon_- m_-,(m_{\tilde P})_{\tilde P}),
\rho(-\epsilon_+ m_+, \epsilon_- m_-,(m_{\tilde P})_{\tilde P})\}$
 \\
\qquad\qquad\qquad\qquad\qquad\qquad $\bigoplus_{\tilde\varphi }\tilde\varphi \otimes \St[[2m_{\tilde\varphi}]]$ &
\\
 $1\otimes \St[[4m_+-1]]\oplus \omega_2\otimes \St[[4m_--1]]
$ 
&
$\{\rho(\epsilon_+ m_+, \epsilon_- m_-,(m_{\tilde P})_{\tilde P}),
\rho(\epsilon_+ m_+, -\epsilon_- m_-,(m_{\tilde P})_{\tilde P})\}$
 \\
\qquad\qquad\qquad\qquad\qquad\qquad $\bigoplus_{\tilde\varphi }\tilde\varphi \otimes \St[[2m_{\tilde\varphi}]]$ &
\\
 $\omega_0\otimes \St[[4m_+-1]]\oplus \omega_2\otimes \St[[4m_--1]]
$ 
&
$\{\rho(\epsilon_+ m_+, \epsilon_- m_-,(m_{\tilde P})_{\tilde P}),
\rho(-\epsilon_+ m_+, -\epsilon_- m_-,(m_{\tilde P})_{\tilde P})\}$
 \\
\qquad\qquad\qquad\qquad\qquad\qquad $\bigoplus_{\tilde\varphi }\tilde\varphi \otimes \St[[2m_{\tilde\varphi}]]$ &
\\
 \hline
\end{tabular}
\end{center}

In the above two tables, each $\tilde\varphi$ is the parameter, with symplectic image and degree $>1$, corresponding to a self-dual polynomial $\tilde P\neq X\pm 1$, and $m_{\tilde\varphi}= m_{\tilde P}$.
\label{intro main theorem symplectic and orthogonal}

\item For $G=$ unramified $\mathrm U(2n)$, we have 
\begin{equation*}
\begin{split}
&\text{the packet of 
$\varphi = \bigoplus_{\tilde\varphi }\tilde\varphi \otimes \St[[2m_{\tilde\varphi}]]$ contains the unique}
\\
&
\text{supercuspidal representation $\rho=\rho_y\times \rho_z$, where $\rho_y = \rho_z = \rho((m_{\tilde P})_{\tilde P})$,}
\end{split}
\end{equation*}
where the correspondence $\tilde\varphi\leftrightarrow \tilde P$ is the same as in (i), but now $\tilde\varphi$, as well as $\tilde P$, are allowed to be of degree 1.
\label{intro main theorem unitary}

\end{enumerate}
As a consequence of (\ref{intro main theorem unitary}), if 
$\tilde\pi$ is a depth zero cuspidal representation of $\tilde G = \GL(m,F)$, compactly induced from a self-dual cuspidal type $\tilde\rho$ with $\tilde{\rho}(\varpi)=1$, then $\tilde\pi$ is conjugate-orthogonal (resp. conjugate-symplectic) if $m$ is odd (resp. even).
\qed\end{thm}


We remark that there are similar results, for $G$ being even orthogonal but both $\mathsf {G}_w$ being odd orthogonal in Proposition \ref{main-result-orthogonal-odd-odd}, and for the connected even orthogonal group $G^0$ in Proposition \ref{main-result-connectd-orthogonal}. We do not summarize them here to avoid introducing further notations.

We explain our methodology by first recalling the foundational knowledge summarized in \cite{Lust-Stevens}. First of all, M{\oe}glin's theory on Jordan blocks \cite{Moeg-exhaustion, Moeg-base-change} (as well as \cite{Moeg-endosc-L-param, Moeglin-Renard-non-quasi-split} for non-quasi-split classical groups) relates the number $a$ in the component $\tilde\varphi \otimes \St[[a]]$ of $\varphi$ with the unique two points $\pm s\in \mathbb R$ where the parabolically induced representation of $\tilde\pi_{\tilde\varphi}|\det|^s\times \pi$, from $M = \GL(m,F)\times G$ to a classical group containing $M$ as a Levi subgroup, is reducible. To study the reducibility and locate these points, we apply Bushnell-Kutzko's theory of types \cite{BK-cover} to transfer the reducibility of parabolically induced representations into the reducibility of modules over the Hecke algebra for the cuspidal type $\tilde\rho\times \rho$. The structure of this Hecke algebra is computed in \cite{stevens-supercusp,Stevens-Miya}: it is a generic Hecke algebra on an infinite dihedral group. We can therefore reduce to consider two Hecke algebras for the types $\tilde \rho\times \rho_w$, where $w\in \{y,z\}$, defined respectively on two reductive quotients of $\tilde {\mathcal J}\times \mathcal J$, each with a generator satisfying an explicit quadratic equation \cite{Lusztig-chevalley}, whose linear coefficient is determined up to a sign. These two equations determine the real-parts of the points $s\in \mathbb C$ at which the above parabolically induced representation is reducible \cite{Blondel-Weil}. For now we already know a lot of information about the parameter of a given supercuspidal representation, i.e., their inertial Jordan blocks \cite{BHS}, or equivalently the tuples of numbers $\{a,b\}$, $\{c,d\}$, and $\{a_{\tilde\varphi}, a_{\tilde\varphi'}\}$ attached to $\varphi$ and associated with $\{1,\omega_0\}$, $\{\omega_1,\omega_2\}$, and $\{{\tilde\varphi}, {\tilde\varphi'}\}$ respectively.

The missing information is hence the exact elementwise matching of the above respective tuples of numbers. As one might guess, this information comes from the sign ambiguities of the linear coefficients of the quadratic equations satisfied by the two Hecke algebra generators, and it is indeed the product of these two signs that matters. With the methodology in this article, this product sign $\nu$ actually comes from the difference of two `canonical' normalizations of the intertwining operator between $\tilde{\rho}\times\rho$ and ${}^\sigma\tilde{\rho}\times\rho$. Here the canonicity means that the normalizations depend completely on $\tilde{\rho}$ and $\rho$. More precisely, there are two suitable normalizations, one for each $w\in \{y,z\}$, of the intertwining operator $\tilde\rho\rightarrow {}^\sigma\tilde\rho$ that each can be expressed as some sort of an orthogonality relation between $\tilde{\rho}$ and $\rho_w$. While exact values of their characters surely determine $\nu$, the typically almost symmetry condition on our representations exhibits a simple yet interesting situation: it allows us to play around with the properties of their characters and obtain $\nu$ by checking whether the difference between the orthogonality relations of $\tilde{\rho}\times\rho_y$ and of $\tilde{\rho}\times\rho_z$ is zero or not, without labouring with calculating their character values. See Section \ref{subsection Covering types for classical groups} - \ref{subsection Proof of the main result} for the complete detail of this methodology, especially the key relation (\ref{Final comparison of eigenvalues}) on comparing eigenvalues of Hecke algebra elements to obtain the difference of the two orthogonality relations in (\ref{cancellation in traced_difference_Dy_Dz}).

As one might have already observed, the typically almost symmetry condition on a given representation $\pi$ also restricts its parameter to specific shapes. In particular, one of the numbers in each tuple $\{a,b\}$, $\{c,d\}$, and $\{a_{\tilde\varphi}, a_{\tilde\varphi'}\}$ must be $-1$ (not 0 because of the parity condition). This conclusion is obtained from the reducibility results stated in Proposition \ref{main results on reducibility points}. Roughly speaking, the results assert that, in almost all cases that we concern, among the four points of reducibility contained in the domain (\ref{fundamental domain of reducibility}) of interest inside $\mathbb{C}$, the two non-real points must be purely imaginary and indeed equal with multiplicity 2. This assertion is deduced from the exact value of the sign $\nu$ determined under the typically almost symmetry condition, whose proof is detailed at the beginning of Section \ref{subsection The final argument}. The two real points of reducibility, which is already computed in \cite{Lust-Stevens}, determines completely the parameter of $\pi$.

The last Section \ref{subsection Examples of quadratic-unipotent representations} includes an example for $\SP(4)$ which resolves the ambiguity arising from \cite[Example 9.4]{Lust-Stevens} with our method.


Previous attempts towards explicit LLC for supercuspidal representations of classical groups besides \cite{Lust-Stevens} include \cite{BHS} for symplectic groups, \cite{Lust-GSP4} for GSp(4), and \cite{Blondel-Tam-2020} and \cite{tam-very-cusp} for unitary groups in specific cases. Much of the ideas of our present paper, including the usages of covering types and Hecke algebras, are based on these articles. Other attempts using explicit character calculations for supercuspidal representations include \cite{Blasco-u2} for U(1,1) and \cite{Adler-Lansky-unram, Adler-Lansky-ram} for U(2,1). Beyond classical groups, Kaletha established a version of LLC for tamely ramified reductive groups \cite{Kaletha-regular-supercuspidal}; when specified to general linear groups, \cite{Oi-Tokimoto} proved that the underlying data of his LLC is the same as those obtained by the author \cite{tam-thesis} using the essentially tame LLC \cite{BH-ET1,BH-ET2,BH-ET3}. With the theory of \cite{AMS-Springer-corresp-cusp-parameter}, Aubert et al. \cite{Aubert-Tsai-Yu-G2} has (very recently) explicated the LLC for the exception group G2.

\subsection*{Acknowledgements}

This research is funded by Radboud Excellence Initiative and the open competition of NWO under Grant No. OCENW.M20.132. The author would like to thank Corinne Blondel and Shaun Stevens for their comments to improve the paper substantially, and Maarten Solleveld for discussions and encouragements.

\subsection*{Notations}



We denote by $F$ a non-Archimedean local field, with ring of integers $ \mathcal O_F$ and the maximal ideal $\mathcal P_F$ generated by a fixed prime element $\varpi=\varpi_F$. Suppose that $\mathcal O_F/\mathcal P_F$ is the finite field $\mathbb F_q$ of $q$ elements, a power of an odd prime number $p$. We usually identify $\mathbb F_q^\times$ with the subgroup $\mu_F$ of the $(q-1)$th roots of unity in $F$.

We use the notation $|\cdot|$ to denote the cardinality of a set, the Archimedean norm of a complex number, or the normalized p-adic norm of a number in $F$ such that $|\varpi|=q^{-1}$; which one being used will be clear from the context. If $x\in \mathbb{C}^\times$, we denote $\sgn(x) = x/|x|$.

Characters of $F^\times$ will be identified with characters of the Weil group $\mathcal W_F$ implicitly via local class field theory. With a fixed $N\in \mathbb N$, the notation $|\det|$ stands for the character $\GL(N,F)\rightarrow \mathbb C^\times$, $x\mapsto |\det x|$. The contragredient of a smooth representation $\pi$ of an $F$-reductive group or the Weil group $\mathcal W_F$ is denoted by $\pi^\vee$.

Transpose of a matrix $A$ is denoted by ${}^t A$, and entry-wise Galois action on $A$ is denoted by ${}^gA$ if $A$ has entries in a field $\bar k$ and $g\in \Gal(\bar k/k)$. If $(A_1,\dots,A_n)$ is  a sequence of either scalars or square matrices, we denote by $\diag(A_1,\dots,A_n)$ (resp. $\antidiag(A_1,\dots,A_n)$) the diagonal (resp. anti-diagonal) matrix with entries $(A_1,\dots,A_n)$ on the diagonal (resp. anti-diagonal, from upper-right to lower-left).

We use bold fonts $\boldsymbol{\mathsf G}$ and $\mathbf G$ to stand for algebraic groups over a finite field and a local field respectively, and normal fonts $\mathsf G$ and $G$ to stand for their subgroups of rational points. When we write $\mathsf G = \GL(N)$ and  $ G = \GL(N)$ for example, we implicitly mean $\mathsf G = \GL(N,\mathbb F_q)$ and  $ G = \GL(N, F)$. Moreover, unless otherwise specified, when we use the term `cuspidal' or `supercuspidal', the described representation is presumed to be irreducible.

\section{Connected reductive groups}

We first provide the preliminary knowledge on cuspidal and supercuspidal representations of general reductive groups (Sec \ref{subsection (Super-)cuspidal representations}). After recalling the basic structures of classical groups (Sec \ref{section Classical groups} and \ref{subsection Finite classical groups}). for the groups over finite fields, we provide a classification on their cuspidal representations (Sec \ref{subsection Finite classical groups}); for the groups over non-Archimedean local fields, we recall the constructions of their depth zero supercuspidal representations via inducing cuspidal types (Sec \ref{subsection Compact inductions}), and describe their parameters via the local Langlands correspondence (LLC) (Sec \ref{subsection Langlands Parameters} and \ref{section LLC for classical groups}). The last subsection \ref{subsection The main result} contains the main results of this paper.

\subsection{(Super-)cuspidal representations}
\label{subsection (Super-)cuspidal representations}

Let $\boldsymbol{\mathsf{G}}$  be a reductive group over the finite field $\mathbb{F}_q$. We put $\mathsf G = \boldsymbol{\mathsf{G}}(\mathbb{F}_q)$ the group of $\mathbb{F}_q$-points, and use similar notations for other reductive groups. Let $\boldsymbol{\mathsf{G}}^0$ be the connected component of $\boldsymbol{\mathsf{G}}$. We denote by $\hat{\boldsymbol{\mathsf{G}}}^0$ the dual group of $\boldsymbol{\mathsf{G}}^0$, the $\mathbb{F}_q$-reductive group whose basic root datum is the dual of that of $\boldsymbol{\mathsf{G}}^0$.

Let $\boldsymbol{\mathsf{T}}$ be a maximal torus in $\boldsymbol{\mathsf{G}}^0$ defined over $\mathbb{F}_q$, and let $\hat{\boldsymbol{\mathsf{T}}}$ be a torus in $\hat{\boldsymbol{\mathsf{G}}}^0$ dual to $\boldsymbol{\mathsf{T}}$. An element $s\in \hat{ {\mathsf{T}}}$ then corresponds to a character $\theta$ of $\mathsf T$, giving a bijection between $\mathsf G^0$-conjugacy classes of pairs $(\mathsf T,\theta)$ with $\hat {\mathsf G}^0$-conjugacy classes of pairs $(\hat{ \mathsf T}, s)$. 

The space of class functions on $\mathsf{G}^0$ is equipped with a Hermitian product $\left<\cdot,\cdot\right>$ such that the basis $\mathcal E(\mathsf{G}^0)$ consisting of isomorphism classes of irreducible (complex) representations of $\mathsf{G}^0$ is orthonormal. Denote by 
$R_{\mathsf{T}}(\theta)$ the Deligne-Lusztig (virtual) representation of $\mathsf G^0$ associated with a character $({\mathsf{T}},\theta)$ corresponding to $(\hat {\mathsf{T}}, s)$ for some $s\in \hat {\mathsf{T}}$. If we define $\mathcal E({\mathsf{G}}^0,s)$ to be the subset  of irreducible representations $\pi\in \mathcal E({\mathsf{G}}^0)$ such that $\left<\pi,R_{\mathsf{T}}(\theta)\right>\neq 0$ for some $({\mathsf{T}},\theta)$, then there is a partition
$$\mathcal E(\mathsf{G}^0) = \bigsqcup_s \mathcal E(\mathsf{G}^0,s)$$
where $s$ runs through the set of $\hat {\mathsf{G}}^0$-conjugacy classes of semisimple elements in $\hat {\mathsf{G}}^0$. We call $\pi\in \mathcal E({\mathsf{G}}^0)$ cuspidal if $\left<\pi,R_{\mathsf{T}}(\theta)\right>=0$ for all pairs $({\mathsf{T}},\theta)$ with $\boldsymbol {\mathsf{T}}$ contained in a proper parabolic subgroup of $\boldsymbol{\mathsf{G}}^0$ defined over $\mathbb{F}_q$.

We sometimes consider representations of a disconnected group. Assume that ${\boldsymbol{\mathsf{G}}}/{\boldsymbol{\mathsf{G}}}^0$ is abelian, and put $\hat{\boldsymbol{\mathsf{G}}} := \hat{\boldsymbol{\mathsf{G}}}^0 \rtimes ({\boldsymbol{\mathsf{G}}}/{\boldsymbol{\mathsf{G}}}^0)^\wedge$, see \cite[Sec 1.3]{Waldspurger-Lusztig-conj} for the precise action of the Pontryagin dual $({\boldsymbol{\mathsf{G}}}/{\boldsymbol{\mathsf{G}}}^0)^\wedge$ of ${\boldsymbol{\mathsf{G}}}/{\boldsymbol{\mathsf{G}}}^0$ on $\hat{\boldsymbol{\mathsf{G}}}^0 $. We then define $\mathcal E(\mathsf{G},s)$, where $s\in \hat{\mathsf{G}}^0$ is semisimple, to be those $\pi\in \mathcal E(\mathsf{G})$ such that there exists an (hence all) irreducible component $\pi^0$ in $\pi|_{\mathsf G^0}$ with $\pi\in \mathcal E(\mathsf{G}^0,s')$ and $s'$ is $\hat {\mathsf{G}}$-conjugate to $s$. Then 
$$\mathcal E(\mathsf{G}) = \bigsqcup_s \mathcal E(\mathsf{G},s)$$
where $s$ runs through the set of $\hat {\mathsf{G}}$-conjugacy classes of semisimple elements in $\hat {\mathsf{G}}^0$.

Now let $\bf G$ be a connected reductive group over a non-Archimedean local field $F$ whose residual field is $\mathbb{F}_q$, and denote by $G$ its group of $F$-points. We will mostly concern about irreducible supercuspidal representations of $G$ of depth zero, whose construction is given in \cite{Moy-Prasad-Jacquet} as follows. Let $\mathcal J_0$ be a maximal parahoric subgroup of $G$, and $\rho_0$ be an irreducible cuspidal representation of $\mathcal J_0$. Let $\mathcal J$ be the normalizer of $\mathcal J_0$ in $G$, and $\rho$ be an irreducible representation of $\mathcal J$ whose restriction to $\mathcal J_0$ contains $\rho_0$. \cite[Prop 6.8]{Moy-Prasad-Jacquet} asserts that the compact induction $\cInd_{\mathcal J}^G \rho$ is an irreducible depth zero supercuspidal representation, and any such representation arises in this way.

\subsection{Classical groups}
\label{section Classical groups}

 Given a finite dimensional vector space $V$ over $F$, we denote by $\tilde{G} = \tilde{G}_V$  the group of $F$-linear automorphisms of $V$. We usually identify $V$ with $F^{\oplus N}$ with $N=\dim _FV$ by choosing a basis, then $\tilde{G}$ is identified with the general linear group $\GL(N)$. The choice of such a basis is usually unimportant and will be made implicitly.

For the setup of classical groups, we require some notations for quadratic extensions. For an extension $F/\Fo$ of degree $\leq 2$, denote by $c$ the Galois involution, taken to be non-trivial when $[F:\Fo]=2$. If further $F/\Fo$ is unramified, then we take a normalizer $\varpi\in \Fo$. The case that $F/\Fo$ is quadratic ramified is not considered in this paper.

Suppose now that $V$ is equipped with a non-degenerate $\epsilon$-Hermitian form $h_V$, relative to an extension $F/\Fo$ of degree $\leq 2$ and with $\epsilon=\pm 1$, which means that
$h_V(w,v)=\epsilon\cdot {}^ch_V(v,w)$ for all $v,w\in V$. We denote by $G$ the subgroup  $\tilde{G}$ of $h_V$-isometries:
$${g\in G}\quad\Leftrightarrow \quad
h_V(gv,gw)= h_V(v,w)\quad \text{ for all }v,w\in V.$$
Hence $G$ is either orthogonal ($\epsilon=1$ and $F=\Fo$), symplectic ($\epsilon = -1$  and $F=\Fo$), or unitary ($\epsilon = (-1)^{N-1}$  and $[F:\Fo]=2$). If $G$ is orthogonal, denote by $G^0$ the subgroup of $g\in G$ with $\det g=1$.

Denote by $\sigma $ the involutive automorphism on $\tilde{G}$ such that  $\tilde{{G}}^\sigma = {{G}}$, so that if we identify $\tilde{G}$ with $\GL(N,F)$ and $h_V$ with $h_J:(v,w)\mapsto {}^{tc}vJw$, where $J\in \tilde{G}$ such that $J = \epsilon{}^{tc}J$, then ${}^\sigma g = J^{-1}{}^{tc}g^{-1}J$ and $G$ consists of $g\in \GL(N)$ such that $g = {}^\sigma g $. 


Let $H$ be the orthogonal sum of $m$ pieces of the 2-dimensional hyperbolic $F$-space of the same type (orthogonal, symplectic, or unitary) as $(V,h)$, i.e., $H=H_-\oplus H_+$, where $H_-= F^{\oplus m}$, and $H_+$ be the $F/\Fo$-Hermitian linear dual of $H_-$, such that $H$ is equipped with the form
\begin{equation*}
h_H((z_-,z_+),(w_-,w_+)) =  (z_-,w_+) +\epsilon{}(w_-,z_+),
\quad 
z_-,w_-\in H_-\text{ and }z_+,w_+\in H_+.
\end{equation*} 
If we identify $H_+$ also with $F^{\oplus m}$ such that 
$(z_-,z_+) = {}^{tc}z_-Kz_+$ for some $K\in \GL(m,F)$, then $h_H$ is given by 
$$(z_-,z_+),(w_-,w_+) \mapsto {}^{tc}z_-Kw_+ +\epsilon{}^{tc}w_-{}^{tc}K z_+.$$ 
We again denote the involutive morphism by $\sigma:\tilde G_{H_-}\rightarrow \tilde G_{H_+}$. The orthogonal sum $W := V \perp H$ is hence equipped with the form $h_W = h_V \perp h_{H_-}$, with the involutive automorphism on $\tilde{G}_{W}$ again denoted by $\sigma$. The corresponding isometry group $G_W$ is of the same type as $G$.

Fix the chain of subspaces $H_-\subset V \oplus H_- \subset W$, and denote by $P$ the corresponding parabolic subgroup in $G_W$ with unipotent radical $N$, and their opposites $P_-$ and $N_-$. Denote by $M$ the Levi component of $P$ that stabilizes the subspaces in the chain, and fix an embedding 
$$i_M:\tilde G_{H_-}\times G_V\rightarrow M, \qquad (g,h)\mapsto \diag(g,h,^\sigma g).$$ 
Denote an element of the unipotent radical $N$ of $P$ by $ (X,Y,Z)^+:=\left[\begin{smallmatrix}
I&X&Y \\ &I&Z \\ &&I
\end{smallmatrix}\right]$, then the involution $\sigma$ on $W$ induces certain relations between the entries $X,Y$ and $Z$. More precisely, we define the $\alpha$-operators on the entries by the following formulas:
\begin{equation}
\label{alpha_on_X_and_Y}
{}^\alpha X :=  -J^{-1} {}^{tc} XK,\qquad
{}^\alpha Y := -\epsilon  {}^{tc}K^{-1} {}^{tc} Y K,
\qquad\text{and}\qquad
{}^\alpha Z := -\epsilon {}^{tc}K^{-1}   {}^{tc} Z J,
\end{equation}
then an element $(X,Y,Z)^+\in N$ if and only if $Z = {}^\alpha X$ and 
\begin{equation}
\label{XalphaX=Y-alphaY}
X{}^\alpha X=Y-{}^\alpha Y,
\end{equation}
and we simply denote such an element by $(X,Y)^+$. Similarly, an element in  $N_-$ takes the form $(X,Y)^-:=\left[\begin{smallmatrix}
I&& \\ {}^\alpha X &I&\\ Y &X&I
\end{smallmatrix}\right]$, where $X$ and $Y$ satisfies the same relation (\ref{XalphaX=Y-alphaY}). Note that $\alpha\circ \alpha$ is the identity on all entries in (\ref{alpha_on_X_and_Y}).

In our present paper, we only consider $G$ to be either an even orthogonal, a symplectic, or a unramified unitary group.

\subsection{Classical groups over finite fields}
\label{subsection Finite classical groups}

Classical groups over $\mathbb F_q$ are defined similarly as in the previous section, by replacing $F/\Fo$ with $\mathbb F_q/\mathbb F_{\qo}$ where $\qo=q$ when ${\mathsf {G}}$ is an orthogonal group $\mathrm O(2n+1)$, $\mathrm O_\delta(2n)$ (here $\delta \in \mathbb F_q^\times /\mathbb F_q^{\times2}\cong \{\pm 1\}$ is the discriminant mod $\mathbb F_q^{\times2}$), or a symplectic group $\SP(2n)$, and $\qo^2=q$ when ${\mathsf {G}}$ is a unitary group $\mathrm{U}(N)$. We use the same notation $c\in \Gal(\mathbb F_q/\mathbb F_{\qo})$ for the descent of $c\in \Gal(F/\Fo)$.

Building on the notion of Lusztig series, we parametrize cuspidal representations of a general linear group $\tilde{\mathsf G} = \GL(m)$ by monic irreducible polynomials in $\mathbb F_q[X]$ of degree $m$. Namely, for a semi-simple element $\tilde s\in \hat{\tilde{\mathsf G}} = \GL(m)$, we have $\mathcal E(\tilde s)$ contains a (unique) cuspidal representation if and only if the (monic) characteristic polynomial $P_{\tilde s}$ of $\tilde s$ is irreducible.

With the same notion, we parametrize cuspidal representations of classical groups over $\mathbb F_q$ by monic polynomials in $\mathbb F_q[X]$ satisfying certain self-duality conditions. The treatment for orthogonal groups will follow the idea from \cite[Ch 4]{Waldspurger-Lusztig-conj} (see also \cite[Lemma 8.9]{lusz-classical} and \cite[Sec 7]{Lust-Stevens}).

If a Lusztig series $ \mathcal E(s)$, where $s\in \hat{\mathsf G}^0$, contains a cuspidal representation of $\mathsf G$, then the (monic) characteristic polynomial of $s$ (when $s$ is viewed as a linear automorphism of the underlying $\mathbb F_q$-space $\hat{\mathsf V}$ of $\hat{\mathsf G}$) factors as 
$P_s = \prod_{\tilde P}\tilde P^{N_{\tilde P}}$, where ${\tilde P}$ ranges over all  monic irreducible polynomials in $\mathbb F_q$ such that 
\begin{enumerate}[(i)]

\item each ${\tilde P}$ has to be self-dual, in the sense that if $\tilde s$ is a root of ${\tilde P}$, then so is ${}^c \tilde s ^{-1}$; 

\item each power $N_{\tilde P}$ takes the following forms:
\begin{enumerate}[(I)]

\item if $\tilde s\neq \pm 1$ or when $\mathsf G$ is unitary, then $N_{\tilde P} = (m_{\tilde P}^2+m_{\tilde P})/2$ for some $m_{\tilde P}\in \mathbb N_{\geq 0}$;

\item If $\tilde s=1$, then $N_{\tilde P} = 2(m_+^2+m_+)+(1-\epsilon)/2$ for some $m_+\in \mathbb N_{\geq 0}$ when ${\mathsf{G}}$ is odd orthogonal or symplectic, and $N_1 = 2m_{+}^2$ for some $m_+\in \mathbb Z$ when ${\mathsf{G}}$ is even orthogonal (of type $(-1)^{m_+}$); 

\item if $\tilde s=-1$, then  $N_{\tilde P} = 2(m_-^2+m_-)$ for some $m_-\in \mathbb N_{\geq 0}$ when ${\mathsf{G}}$ is odd orthogonal, and $N_{-1} = 2m_-^2$ , for some $m_-\in \mathbb Z$ when ${\mathsf{G}}$ is symplectic or even orthogonal (of type $(-1)^{m_-}$).
\end{enumerate}\end{enumerate}
When $q=\qo$, we denote by  $\rho = \rho(m_+,m_-,(m_{\tilde P})_{\tilde P})$ a cuspidal representation corresponding to $P$ as above, where ${\tilde P}$ ranges over self-dual irreducible  polynomials over $\mathbb F_q$ of degree $>1$. There are either 1, 2, or 4 such representations with a fixed $P$, depending on the type of $\mathsf G$ and whether $m_+$ and $m_-$ are zero or not. Here $\rho$ is a representation of $\mathsf G^0$ when $\mathsf G$ is odd orthogonal, and is of $\mathsf G$ when $\mathsf G$ is even orthogonal. Note that the numbers in $(m_+,m_-,(m_{\tilde P})_{\tilde P})$ determines the dimension of the Hermitian space $\mathsf{V}$ defining $\mathsf{G}$, and the pair $(m_+,m_-)$ determines the parity $\delta = (-1)^{m_++m_-}$ when $\mathsf{G}$ is even orthogonal.

When $\mathsf G$ is unitary, we use similar notations $\rho = \rho((m_{\tilde P})_{\tilde P})$ to denote the cuspidal representation corresponding to $P$, where ${\tilde P}$ ranges over self-dual irreducible  polynomials over $\mathbb F_q$. There is a unique such representation with a fixed $P$.

When  $\mathsf{G} = \mathsf{G(V)}$  is even orthogonal or symplectic, we provide a few properties of these representations which will be used later on. Denote by $\omega$ the character of  $\mathbb{F}_q^\times$ of order 2, and by $\mathsf{GG}(\mathsf V)$ the group of (orthogonal or symplectic, depending on $\mathsf{G}$) similitudes, with a factor map denoted by $\lambda:\mathsf{GG}(\mathsf V)\rightarrow \mathbb{F}_q^\times$. Take $h\in \mathsf{GG}(\mathsf V)$ such that $\omega\circ \lambda(h)=-1$.

\begin{lem}
 \label{properties as in Waldspurger book} 
Fix non-negative integers $m_{\tilde P}$, with $\tilde P \in \mathbb F_q[X]$ self-dual of degree $>1$. Let $(m_+,m_-)\in \mathbb{N}_{\geq 0}\times \mathbb{Z}$ in the symplectic case, and $(m_+,m_-)\in \mathbb{Z}\times \mathbb{Z}$ in the even orthogonal case.
\begin{enumerate}[(i)]

\item We have 
$$\rho(m_+,m_-,(m_{\tilde P})_{\tilde P})\circ \Ad(h) \cong \rho(m_+,-m_-,(m_{\tilde P})_{\tilde P}),$$
 and if $\mathsf{G}$ is even orthogonal, then we have furthermore
 $$\rho(m_+,m_-,(m_{\tilde P})_{\tilde P})\otimes \det \cong \rho(-m_+,-m_-,(m_{\tilde P})_{\tilde P}).$$ 
 Note: when all $m_{\tilde P}=0$ for $\tilde P\neq X\pm 1$, these two properties are proved in \cite[4.3]{Waldspurger-Lusztig-conj}. 
 \label{concerning det and conjugate by h}
 
\item With $(|m_+|,|m_-|)$ also fixed, all $\rho( m_+, m_-,(m_{\tilde P})_{\tilde P})$ have the same dimension, i.e., the dimension is independent of the signs of $(m_+, m_-)$.
 \label{concerning dimension}
\end{enumerate}
\end{lem}
\proof
We first recall the lifting map from \cite[7.8 and 7.9]{lusz-classical}, for $s\in \hat {\mathsf G}^0$ the dual group of $\mathsf G^0$, 
$$\mathcal R_s:\mathcal E (Z_{\hat {\mathsf G}}(s)^\wedge,1) \rightarrow \mathcal E ({\mathsf L},s)\rightarrow \mathcal E ({\mathsf G},s).$$
Here ${\mathsf L}$ is taken to be the dual group of 
$$\hat {\mathsf L}:=\mathsf G(\hat {\mathsf V}_{\pm 1}) \times \prod_{{\tilde s}\neq \pm 1}\mathsf G(\hat {\mathsf V}_{\tilde s})\hookrightarrow \hat {\mathsf G} = \hat {\mathsf G}(\hat {\mathsf V}),$$
where $\hat{\mathsf V}_{\tilde s}$ is the $\tilde s$-eigen-subspace of $\hat{\mathsf V}$ for all eigenvalues $\tilde s$ appeared in $s$, so that $\mathsf G(\hat {\mathsf V}_{\pm 1})$, where $\hat{\mathsf V}_{\pm 1} = \hat{\mathsf V}_{1}\oplus \hat{\mathsf V}_{- 1}$, is a classical group of the same type as $\hat {\mathsf G}$, and $G(\hat {\mathsf V}_{\tilde s})$ are unitary groups for all $\tilde s\neq \pm 1$. 

We now describe the two maps composing $\mathcal R_s$. First of all, these two maps are bijections. They map cuspidal representations to cuspidal representations, and also conversely. For the first map $\mathcal E (Z_{\hat {\mathsf G}}(s)^\wedge,1) \rightarrow \mathcal E ({\mathsf L},s)$, since 
$$Z_{\hat {\mathsf G}}(s) = Z_{\hat {\mathsf L}}(s)  = \mathsf G(\hat {\mathsf V}_{ 1})\times \mathsf G(\hat {\mathsf V}_{- 1}) \times \prod_{{\tilde s}\neq \pm 1}\mathsf G(\hat {\mathsf V}_{\tilde s}),$$
it suffices to assume that $s^2=1$, i.e., we are reduced to consider cuspidal quadratic unipotent representations. The map can then be essentially defined by \cite[Lemma 7.9]{lusz-classical}, and is described more precisely in \cite[Proof of Prop 4.3]{Waldspurger-Lusztig-conj}. Briefly speaking, when $\mathsf G$ is even orthogonal, the set $\mathcal E (Z_{\hat {\mathsf G}}(s)^\wedge,1)$ consists of four cuspidal representations when both $m_+$ and $m_-$ are non-zero, two when exactly one of $m_+$ or $m_-$ is non-zero, and one if both $m_+$ and $m_-$ are zero. For the simplicity of the writing, we ignore the last two special cases. When both $m_+,m_-\neq 0$, these four cuspidal representations are lifted to the four cuspidal representations in $\mathcal E ({\mathsf L},s)$, which satisfy the two properties in \cite[4.3.Prop(iii),(v)]{Waldspurger-Lusztig-conj}, i.e., the properties in (i) when all $m_{\tilde P}=0$ for $\tilde P\neq X\pm 1$. Also, we can easily check the dimension formulas in \cite[9.9.Cor(i), 13.9.Prop(i)]{Waldspurger-Lusztig-conj} that the dimensions are independent of the signs, i.e., they only depend on $|m_+|$ and $|m_-|$ (or more precisely, in the language of \emph{loc. cit.}, dimensions are independent under the exchange of $\Lambda_i^+$ and $\Lambda_i^-$ for a given pair of symbols $(\Lambda_i)_{i=1,2}$).

The second map is, up to a sign, the Deligne-Lusztig induction $ \pm R_{\mathsf L}^{\mathsf G}:\mathcal E ({\mathsf L},s)\rightarrow \mathcal E ({\mathsf G},s)$. (The disconnected analog is defined in \cite{DM-non-conn}.) Since $\mathsf G/\mathsf G^0\cong \mathsf L/\mathsf L^0 \cong \mathsf G(\hat {\mathsf V}_{\pm 1})/\mathsf G(\hat {\mathsf V}_{\pm 1})^0$, we can take an element $h_{\pm 1}$ in the group $\mathsf G\mathsf G_{\pm 1}$ of similitudes of $\mathsf G_{\pm 1} = \mathsf G(\hat {\mathsf V}_{\pm 1})$, defined analogous as the element $h\in \mathsf{GG} = \mathsf{GG(V)}$, such that  $h_{\pm 1}\mapsto h$ via the canonical map $\mathsf G\mathsf G_{\pm 1}\rightarrow \mathsf{GG}$. The map $ \pm R_{\mathsf L}^{\mathsf G}$ then lifts the four cuspidal representations in $\mathcal E ({\mathsf L},s)$ to the four cuspidals in $\mathcal E ({\mathsf G},s)$ satisfying the same properties of our proposition, i.e., 
\begin{equation*}
\begin{split}
(\pm R_{\mathsf L}^{\mathsf G}\rho)\otimes \det{}_{\mathsf G}&\cong 
\pm R_{\mathsf L}^{\mathsf G}(\rho \otimes \det{}_{\mathsf G(\hat {\mathsf V}_{\pm 1})}),
\\
(\pm R_{\mathsf L}^{\mathsf G}\rho)\circ \Ad(h)
&\cong 
\pm R_{\mathsf L}^{\mathsf G}(\rho \circ \Ad (h_{\pm 1})).
\end{split}
\end{equation*}
This completes the proof of (i) when $\mathsf G$ is even orthogonal. The arguments of (i) for $\mathsf G$ symplectic is similar, by replacing the number of cuspidal representations from four to two. The proof of (ii) follows from the formula in the proof of \cite[Lemma 7.9]{lusz-classical}: 
$\dim R_{\mathsf L}^{\mathsf G}\rho = |\mathsf G|_{p'}|\mathsf L|_{p'}^{-1}\dim \rho$, 
and the fact that $\mathsf L$ depends only on  $|m_+|$ and $|m_-|$, but not on the signs.
\qed

We will need a parametrization of cuspidal representations of $\mathsf G$ when it is odd orthogonal. Denote $\rho = \rho(m_{+},m_{-}, (m_{\tilde P})_{\tilde P})$, with all $m_{+},m_{-}$ and $ m_{\tilde P}$ non-negative. Then $\rho$ extends in two ways from $\mathsf{G}^0$ to $\mathsf{G}= \mathsf{G}^0\times\{\pm1\}$, denoted by $\rho(+)$ and $\rho(-)$ in the obvious way.

When $\mathsf G$ is even orthogonal, we will also need a parametrization of cuspidal representations of $\mathsf G^0 $ via that of $\mathsf G $. Given $m_+,m_-\in \mathbb Z$, Lemma \ref{properties as in Waldspurger book}(\ref{concerning det and conjugate by h}) implies that, when $(m_+,m_-)\neq (0,0)$, the restriction of  $\rho(m_+,m_-,(m_{\tilde P})_{\tilde P}) $ to $\mathsf G^0$ is equal to that of $\rho(-m_+,-m_-,(m_{\tilde P})_{\tilde P}) $, and so is irreducible. We denote this restriction by 
$\rho(\pm (m_+,m_-), (m_{\tilde P})_{\tilde P})$. When $(m_+,m_-)=(0,0)$, the restriction of   $ \rho(0,0,(m_{\tilde P})_{\tilde P}) $ splits into two irreducible representations, and we label them by $\rho_+(0,0,(m_{\tilde P})_{\tilde P})$ and $\rho_-(0,0,(m_{\tilde P})_{\tilde P})$.

\subsection{Compact inductions}
\label{subsection Compact inductions}

We now describe the constructions of depth zero supercuspidal representations of general linear groups and classical groups over $\Fo$ via compact inductions of cuspidal types.

A maximal compact subgroup $\tilde{\mathcal J}_0$ of $\tilde G = \GL(m)$ is conjugate to $\GL(m,\mathcal O_F)$. Let $\tilde{\mathcal J}_{0+}$ be the pro-p-radical of $\tilde{\mathcal J}_0$, so that $\tilde{\mathcal J}_{0}/\tilde{\mathcal J}_{0+}\cong \tilde{\mathsf G}:= \mathsf{GL}(m)$. Given $\tilde\rho$ a cuspidal representation of $\tilde{\mathsf G}$, we inflate $\tilde{\rho}$ to $\tilde{\mathcal J}_0$ and extend it to a representation  on the normalizer $\tilde{\mathcal J}:=F^\times \tilde{\mathcal J}_0$ of $\tilde{\mathcal J}_0$ in $\tilde G$, by multiplying an extended character $\tilde{\xi}$ from $F^\times$ to  $\tilde{\mathcal J}$ agreeing with $\tilde\rho$ on $\tilde{\mathcal J}_0\cap F^\times  = \mathcal O_F^\times$. Denote the extension on $\tilde{\mathcal J}$ by $\tilde{\rho}_{\tilde{\xi}}$, then the compactly induced representation
\begin{equation}
\label{depth zero supercuspidal compact induction}
\tilde \pi_{\tilde{\rho}_{{\tilde{\xi}}
}}:=\cInd_{\tilde{ \mathcal J}}^{\tilde G}\tilde{\rho}_{\tilde{\xi}}
\end{equation}
is a supercuspidal representation of $\tilde G $.

The construction of supercuspidal representations of a classical group $G $ is similar, but the precise description branches into cases. We first recall, from \cite[Sec 3.3]{stevens-supercusp} that the reductive quotient of a maximal compact subgroup $\mathcal J$ in $G$ is a product $\mathsf G_y\times \mathsf G_z$ of two classical groups over $\mathbb F_q$, both of the same type as $G$, and with one of them allowed to be trivial. (Note that we have excluded $G$ to be a ramified unitary group, whose maximal compact subgroups have reductive quotients of the form a product of a symplectic group and an orthogonal group.) When $G \supsetneq G^0$, i.e., when $G$ is an orthogonal group, the maximal compact subgroups in $G^0$ are obtained by intersecting those in $G$ with $G^0$. In general, these maximal compact subgroups properly contain their underlying parahoric subgroups.

For constructing supercuspidal representations of classical groups, we require our maximal compact subgroups whose underlying parahorics are also maximal; see the discussion in \cite[Appendix]{Stevens-Miya}. This requirement only excludes the possibility that one of the $\mathsf G_w$, where $w\in \{y,z\}$, is of the type the split SO(2). These maximal compacts subgroups are precisely the normalizers of maximal parahoric subgroups.

For $w\in\{y,z\}$, let $\rho_w$ be a cuspidal representation of $\mathsf G_w$, denote $\rho = \rho_y \times \rho_z$, and inflate $\rho$ to a representation of a maximal compact subgroup $\mathcal J$ whose reductive quotient is $\mathsf G_y\times \mathsf G_z$. The compact induction 
$$ \pi_{{\rho}_{}}:=\cInd_{{ \mathcal J}}^{G}{\rho}_{}$$
is a supercuspidal representation of $G$. Suppose now $G\supsetneq G^0$, which happens when $G$ is orthogonal, $\mathcal J^0:=G^0\cap \mathcal J$, and $\rho^0$ is a component of $\rho|_{\mathcal J^0}$. We have either $\rho^0 = \rho|_{\mathcal J^0}$, or there exists ${\rho^0}'$ such that $ \rho|_{\mathcal J^0} = \rho^0 \oplus {\rho^0}'$. The compact inductions 
$ \pi_{{\rho}^{0}}:=\cInd_{\mathcal J^0}^{G^0}{\rho}^{0}$ in both cases and $\pi_{{{\rho}^{0}}'}:=\cInd_{{ \mathcal J}}^{G^0}{{\rho}^{0}}'$ in the second case
are supercuspidal representations of $G^0$.  Since $\mathcal J/\mathcal J^0\cong G/G^0$, we have  $ \pi_{\rho}|_{G^0}= \pi_{\rho^0}  $ in the first case, and $  \pi_{\rho}|_{G^0}  = \pi_{{\rho}^{0}} \oplus \pi_{{\rho}^{0'}} $ in the second.

We require a few more details on the inducing types that allow us to obtain a complete classification. Let's first consider when both $\mathsf{G}_w$, for $w\in \{y,z\}$, are odd orthogonal. Denote $\rho_w = \rho(m_{+,w},m_{-,w}, (m_{\tilde P,w})_{\tilde P})$, with all $m_{+,w},m_{-,w}$ and $ m_{\tilde P,w}$ non-negative. Each $\rho_w $ extends in two ways $\rho_w(+)$ and $\rho_w(-)$ from $\mathsf{G}^0_w$ to $\mathsf{G}_w = \mathsf{G}^0_w\times\{\pm1\}$. Let $\mathcal J$ be a compact subgroup of $G$ such that $\mathcal J \cap G^0$ contains the parahoric subgroup $\mathcal J^0_0$ in $G^0$ with reductive quotient $\mathcal J^0_0/\mathcal J^0_{0+}\cong \mathsf G^0_y\times\mathsf G^0_z$ and that $\mathcal J/\mathcal J^0_{0+}\cong \mathsf G_y\times\mathsf G_z$. Then $\mathcal J$ is its own normalizer in $G$, and the compact inductions of the inflations of 
\begin{equation*}
\rho_y(\epsilon_y)\times \rho_z(\epsilon_z),\quad \epsilon_y,\epsilon_z\in \{\pm\}
\end{equation*}
 are supercuspidal representations of $G$. Also, the product $\rho_y\times \rho_z$ extends in two ways from $\mathsf{G}^0_y\times \mathsf{G}^0_z$ to
$$\mathsf{S}(\mathsf G_y\times\mathsf G_z):=\{(a,b)\in\mathsf G_y\times\mathsf G_z: \det(a)\det(b)=1 \},$$ 
which is isomorphic to $\mathsf G^0_y\times\mathsf G^0_z\times \{\pm 1\}$. Denoted the extensions, again in the obvious way, by 
\begin{equation}
\label{(++)on-S(OtimesO)-odd-odd}
\rho(\pm):=(\rho_y\times \rho_z)(\pm).
\end{equation}
The compact subgroup  $\mathcal J^0$ in $G^0$ with quotient $\mathcal J^0/\mathcal J^0_{0+}\cong \mathsf{S}(\mathsf G_y\times\mathsf G_z)$ is the normalizer of the maximal parahoric $\mathcal J^0_0$ in $G^0$, and so the compact induction of $(\rho_y\times \rho_z)(\pm)$ to $G^0$ is a supercuspidal representation.

We then consider when both $\mathsf{G}_w$ are even orthogonal. Denote $\rho_w = \rho(m_{+,w},m_{-,w}, (m_{\tilde P,w})_{\tilde P})$ with all $m_{+,w},m_{-,w}\in \mathbb{Z}$ and other $ m_{\tilde P,w}$ non-negative, and form the product type $\rho_y\times \rho_z$ of $\mathsf G_y\times\mathsf G_z $. There are three possibilities.
\begin{enumerate}[(i)]
\item $(m_{+,w},m_{-,w})\neq (0,0)$
for both $w\in \{y,z\}$;

\item $(m_{+,w},m_{-,w})=(0,0)$
for both $w\in \{y,z\}$;

\item exactly one pair of $(m_{+,y},m_{-,y})$ and $(m_{+,z},m_{-,z})$ is $(0,0)$.
\end{enumerate}
We will only need the cases (i) and (ii) in this paper. In the first case when $(m_{+,w},m_{-,w})\neq (0,0)$
for both $w\in \{y,z\}$. Consider the cuspidal representations 
\begin{equation}
\label{(++)(++)on-OtimesO}
\begin{split}
&\rho(m_{+,y},m_{-,y},(m_{\tilde P,y})_{\tilde P})\times \rho(m_{+,z},m_{-,z},(m_{\tilde P,z})_{\tilde P})
\\
\text{ and }\quad&\rho(-m_{+,y},-m_{-,y},(m_{\tilde P,y})_{\tilde P})\times \rho(-m_{+,z},-m_{-,z},(m_{\tilde P,z})_{\tilde P})
\end{split}
\end{equation}
on $\mathsf G_y\times\mathsf G_z$. Let $\mathcal J$ is a compact subgroup of $G$ defined as in the precious case, then the compact inductions of (\ref{(++)(++)on-OtimesO}) from $\mathcal J$ to $G$ are supercuspidal representations.

Both representations in (\ref{(++)(++)on-OtimesO}) restrict to a cuspidal representation on $\mathsf{S}(\mathsf G_y\times\mathsf G_z)$,  
denoted by 
\begin{equation}
\label{(++)on-S(OtimesO)}
\rho(\pm):=\rho(\pm (m_{+,y},m_{-,y};m_{+,z},m_{-,z});(m_{\tilde P,y})_{\tilde P};(m_{\tilde P,z})_{\tilde P}).
\end{equation}
It restricts to the cuspidal representation 
\begin{equation}
\label{type-G0yTtimesG0z}
\text{$\rho(\pm (m_{+,y},m_{-,y}),(m_{\tilde P,y})_{\tilde P})\times \rho(\pm (m_{+,z},m_{-,z}),(m_{\tilde P,z})_{\tilde P})$\quad  on $\mathsf G^0_y\times\mathsf G^0_z$. }
\end{equation}
Let $\mathcal J^0$ be the compact subgroup of $G^0$ defined previously, then the compact induction of $\rho(\pm)$ from ${\mathcal J^0}$ to ${G^0}$ is a supercuspidal representation, the one that extends to $G$ in two ways as the compact inductions of both representations in (\ref{(++)(++)on-OtimesO}).

Similarly, both 
\begin{equation*}
\begin{split}
&\rho(m_{+,y},m_{-,y},(m_{\tilde P,y})_{\tilde P})
\times \rho(-m_{+,z},-m_{-,z},(m_{\tilde P,z})_{\tilde P})
\\
\text{ and }\quad&\rho(-m_{+,y},-m_{-,y},(m_{\tilde P,y})_{\tilde P})\times \rho(m_{+,z},m_{-,z},(m_{\tilde P,z})_{\tilde P})
\end{split}
\end{equation*}
restrict to a cuspidal representation 
\begin{equation}
\label{(+-)on-S(OtimesO)}
\begin{split}
&\rho(\pm (m_{+,y},m_{-,y};-m_{+,z},-m_{-,z});(m_{\tilde P,y})_{\tilde P};(m_{\tilde P,z})_{\tilde P}) 
\\
&= \rho(\pm (-m_{+,y},-m_{-,y};m_{+,z},m_{-,z});(m_{\tilde P,y})_{\tilde P}
;(m_{\tilde P,z})_{\tilde P})
\end{split}
\end{equation}
on $\mathsf{S}(\mathsf G_y\times\mathsf G_z)$, and both (\ref{(++)on-S(OtimesO)}) and (\ref{(+-)on-S(OtimesO)}) restrict to the cuspidal representation (\ref{type-G0yTtimesG0z}).

Suppose that all $(m_{\tilde P,w})_{\tilde P}$, for $w\in\{y,z\}$ and $\deg \tilde P>1$, are fixed. When the four numbers $|m_{\epsilon,w}|$, where $(\epsilon,w)\in \{\pm1\}\times \{y,z\}$, are also fixed, and depending on whether they contain 0,1, or 2 zeros, we have classified 16, 8, or 4 supercuspidal representations of $G$, and 8, 4, or 2 supercuspidal representations of $G^0$.

In the second case when $(m_{+,w},m_{-,w})=(0,0)$
for both $w\in \{y,z\}$, the cuspidal representation 
\begin{equation}
\label{all-zero-on-OtimesO}
\rho(0,0,(m_{\tilde P,y})_{\tilde P})\times \rho(0,0,(m_{\tilde P,z})_{\tilde P})
\end{equation} splits into a direct sum of two cuspidal representations 
under the restriction to $\mathsf{S}(\mathsf G_y\times\mathsf G_z)$, denoted by 
\begin{equation}
\label{all-zero-on-S(OtimesO)}
\rho_\epsilon := \rho_\epsilon(0,0,(m_{\tilde P,y})_{\tilde P}
;0,0,(m_{\tilde P,z})_{\tilde P}
),\quad \epsilon\in \{\pm\}.
\end{equation}
We have classified two supercuspidal representations in $G^0$ by compactly inducing  (\ref{all-zero-on-S(OtimesO)}), and they both induce to the compact induction of the inflation of (\ref{all-zero-on-OtimesO}). Note that since $\rho(0,0,(m_{\tilde P,w})_{\tilde P})$ splits into a direct sum $\rho_+\oplus \rho_-$ upon restriction to $\mathsf G_w^0 $, each (\ref{all-zero-on-S(OtimesO)}) further splits into the direct sum of 
\begin{equation}
\label{all-zero-on-SO-times-SO)}
\rho_{+}\times \rho_{\epsilon}\text{ and }\rho_{-}\times \rho_{-\epsilon}\quad \text{ on }\mathsf G_y^0 \times  \mathsf G_z^0.
\end{equation}

Finally, the case when both $\mathsf{G}_w$ symplectic is similar and simpler. Here $\rho_w = \rho(m_{+,w},m_{-,w}, (m_{\tilde P,w})_{\tilde P})$ with only $m_{-,w}\in \mathbb{Z}$ and other $m_{+,w}$ and $ m_{\tilde P,w}$ non-negative. The compact subgroup $\mathcal J$ with quotient $\mathcal J/\mathcal J_{0+}\cong \mathsf G_y\times\mathsf G_z$ is its own normalizer in $G$, and the compact inductions of $ \rho_y\times \rho_z$
is a supercuspidal representation of $G$.

\subsection{Langlands parameters}
\label{subsection Langlands Parameters}

Let $\mathcal{W}_{\Fo}$ be the Weil group of $\Fo$. For a connected reductive group $G$ defined over $\Fo$, denote by $\hat{\mathbf G}$ its Langlands dual group. A Langlands parameter $\varphi$ for $G = \mathbf{G}(\Fo)$ is a $\hat{\mathbf{G}}$-conjugacy class of a smooth (Frobenius-semisimple) morphism from the Weil-Deligne group ${\mathcal{W}}_{\Fo}\times \SL(2,\mathbb C)$ to the L-group ${}^L \mathbf G:=\hat{\mathbf G} \rtimes {\mathcal{W}}_{\Fo}$ of $\mathbf{ G}$, i.e., $\varphi:{\mathcal{W}}_{\Fo}\times \SL(2,\mathbb C)\rightarrow {}^L \mathbf G$, such that the composition of $\varphi $ with the projection ${}^L \mathbf G \rightarrow {\mathcal{W}}_{\Fo}$ is the identity map on ${\mathcal{W}}_{\Fo}$ and $\varphi|_{\SL(2,\mathbb C)}$ is an algebraic morphism.

For the general linear group $\tilde {G} = \GL(m)$, a parameter can be simply viewed as (an isomorphism class of) a representation of ${\mathcal{W}}_{F}\times \SL(2,\mathbb C)$ of degree $m$. As we only require the supercuspidal part of the LLC for $\tilde G$, we will focus on irreducible representations of ${\mathcal{W}}_{F}$ degree $m$ as parameters for $\tilde G$.

For classical groups, we view a parameter $\varphi$ as a morphism whose image is embedded into $\GL(\hat N,\mathbb C)$ for a suitable $\hat N\in \mathbb N_{\geq 0}$. For the groups we consider in this paper, namely even orthogonal, symplectic, and unramified unitary groups, $\hat N$ is determined according to the following.
\begin{equation*}
\begin{split}
\mathbf G = \mathrm O(2n),\quad &{}^L \mathbf G^0:=\SO(2n,\mathbb C) \rtimes {\mathcal{W}}_{F} \twoheadrightarrow \SO(2n,\mathbb C)\rtimes {\Gal}({E/F})\cong \mathrm O(2n,\mathbb C) \hookrightarrow \GL(2n,\mathbb C);
\\
\mathbf G = \SP(2n),\quad &{}^L \mathbf G:=\SO(2n+1,\mathbb C) \times {\mathcal{W}}_{F} \twoheadrightarrow \SO(2n+1,\mathbb C)\hookrightarrow \GL(2n+1,\mathbb C);
\\
\mathbf G = \mathrm U(N),\quad &\mathrm{image}\varphi\subset  \GL(N,\mathbb C) \times {\mathcal{W}}_{F}  \hookrightarrow {}^L \mathbf G:=\Res_{F/\Fo}\GL(N,\mathbb C) \rtimes {\mathcal{W}}_{\Fo} .
\end{split}
\end{equation*}
For even orthogonal groups, $E/F$ is the quadratic extension determined by the discriminant of the Hermitian form on $V$, so that $\mathcal W_E = \ker \det \varphi$. By factoring through $\mathcal W_F \rightarrow \mathcal W_F / \mathcal W_E \cong \Gal(E/F)\cong \mathrm{O}(2n,\mathbb C)/\mathrm{SO}(2n,\mathbb C)$, we regard the image of $\varphi$ in $\mathrm O(2n,\mathbb C) $ via the isomorphism 
\begin{equation}
\label{SO2nTimesGalIsIsomorphicToO2n}
\SO(2n,\mathbb C)\rtimes {\Gal}({E/F})\cong \mathrm O(2n,\mathbb C). 
\end{equation}
 We remark from \cite[Sec 3.6]{atobe-Gan-SO-even} that this isomorphism is not canonical. The  $\Gal(E/F)$-action on $\SO(2n,\mathbb{C})$ has to induce the order 2 action on the Dynkin diagram: we take 
$$\mathbbm p = \diag(I_{n-1},\antidiag(1,1),I_{n-1})$$
and choose the isomorphism (\ref{SO2nTimesGalIsIsomorphicToO2n}) to be identity on $\SO(2n,\mathbb C)$ and the non-trivial element $\tau\in {\Gal}({E/F})$ is mapped to $\mathbbm p$. Another choice that is also identity on $\SO(2n,\mathbb C)$ is given by $\tau\mapsto -\mathbbm p$.

We usually view $\varphi$ as a finite dimensional representation of ${\mathcal{W}}_{\Fo}\times \SL(2,\mathbb C)$. However, according to \cite[Th 8.1(ii)]{GGP}, the isomorphism class of the representations determines its underlying parameter, or more precisely its $\hat{\mathbf G}^0$-conjugacy class, except when $G$ is even orthogonal and all irreducible components of $\varphi$ are even dimensional. In this exceptional case, the representation determines only the $\hat{\mathbf G}$-conjugacy class of the parameter, and there are two $\hat{\mathbf G}^0$-conjugacy classes of $\varphi$ that give rise to the same representation.

 We denote by $\St[a]$ the $a$-dimensional irreducible representation of $\SL(2,\mathbb C)$. Here $a$ is allowed to be $0$ or $-1$, in which case we put $\St[0] $ and $ \St[-1] $ to be the zero representation, and view them as symplectic and orthogonal representations respectively. Denote 
$$\St[[a]] = \St[a] \oplus \St[a-2]\oplus \cdots.$$
Given an irreducible representation $\tilde{\varphi}$ of ${\mathcal{W}}_F$, we denote shorthand $\tilde\varphi[[a]] = \tilde\varphi\otimes \St[[a]]$.

Parallel to the condition of eigenvalues in Section \ref{subsection Finite classical groups}, we will only consider those parameters with irreducible components of the form $\tilde\varphi[[a]]$ such that $ \tilde\varphi$ is self-dual, which means that $ \tilde\varphi \cong  {}^\sigma \tilde\varphi := {}^c \tilde\varphi^\vee$. Let's first look at self-dual characters, i.e., when $\deg \tilde\varphi =1$.

\begin{enumerate}[(i)]
\item When $F=\Fo$, the self-dual condition on $\tilde\varphi$ means that $\tilde\varphi$ is either trivial or quadratic. We denote by $\omega_0$ the unramified quadratic character of $F^\times$, by $\omega_1$ the ramified quadratic character of $F^\times$ with $\omega_1(\varpi)=1$, and put $\omega_2 = \omega_1\omega_0$. All four characters $1,\omega_0,\omega_1,\omega_2$ are orthogonal (while there is no such thing as a `symplectic character').

\item When $F \neq {\Fo}$, the self-dual condition on $\tilde\varphi$ implies that $\tilde\varphi|_{\Fo^\times}$ is quadratic, and we call $\tilde\varphi$ conjugate-orthogonal (resp. conjugate-symplectic) if $\tilde\varphi|_{\Fo^\times}$ is trivial (resp. equal to $\omega_{F/\Fo}$, the quadratic character corresponding to the extension $F/\Fo$ via local class field theory). Therefore, when $F/\Fo$ is unramified, a self-dual $\tilde\varphi$ is conjugate-orthogonal (resp. conjugate-symplectic) if and only if $\tilde\varphi(\varpi)=1$ (resp. $-1$). 

\end{enumerate}
We then consider an irreducible depth zero  representation  $\tilde\varphi$ with degree $m>1$. Here by depth zero it means that $\tilde\varphi$ is trivial on the wild-inertia subgroup of $\mathcal W_F$. It is known that $\tilde\varphi$ is of the form  $\tilde\varphi \cong \Ind^{\mathcal W_F}_{\mathcal W_E}\tilde\xi$ for some depth zero character $\tilde\xi$ of $E^\times$, where $E/F$ is the unramified extension of degree $m:=\deg\tilde\varphi$. Here $\tilde\xi$ is necessarily regular, i.e., ${}^\sigma\tilde\xi\neq \tilde\xi$ for all non-trivial $\sigma\in \Gal(E/F)$, because $\tilde\varphi$ is irreducible.

If $\tilde\varphi$ is furthermore self-dual, then so is $\tilde\xi$, which means the following. 
\begin{enumerate}[(i)]
\item  When $F=\Fo$, then $m$ is even and we denote $\Eo$ the unramified extension of $F$ of degree $m/2$. Then $\tilde\xi$ being self-dual means that $\tilde\xi = {}^{c}\tilde\xi^{-1}$, where $c\in \Gal(E/\Eo)$ is the non-trivial automorphism. 

\item When $F/\Fo$ is quadratic unramified, then $m$ is odd. We denote $\Eo/\Fo$ the unramified extension of degree $m$, and $\tilde\xi$ being self-dual means the same as in the previous case: $\tilde\xi = {}^{c}\tilde\xi^{-1}$, where $1\neq c\in \Gal(E/\Eo)$.

\end{enumerate}
Note that in both cases, if $\tilde\xi$ is  conjugate-orthogonal (resp. conjugate-symplectic), then $\tilde\varphi$ is orthogonal (resp. symplectic) when $F=\Fo$, and is conjugate-orthogonal (resp. conjugate-symplectic) when $F\neq \Fo$.

Since we focus only on depth zero supercuspidal representations of $G$, we take representations of ${\mathcal{W}}_{F}\times \SL(2,\mathbb C)$ of the following form as our parameters. If $G$ is orthogonal or symplectic, we take
\begin{equation}
\label{quadratic-parameter}
\varphi = 1[[a]]\oplus \omega_0[[b]]\oplus\omega_1[[c]]\oplus \omega_2[[d]]\oplus \bigoplus_{
\tilde\varphi}
\tilde\varphi[[a_{\tilde\varphi}]]\oplus \tilde\varphi'[[b_{\tilde\varphi}]],
\end{equation}
with $\tilde\varphi $ in the sum ranging over  
\begin{equation}
\label{self-dual-condition-deg>1}
\begin{split}
&\text{irreducible representations $ \Ind^{\mathcal W_F}_{\mathcal W_E}\tilde\xi$ such that}
\\
&\text{$E/F$ is unramified of degree $>1$, and }
\\
&\text{$\tilde\xi$ is self-dual, depth zero, and with $\tilde\xi(\varpi)=1$,}
\end{split}
\end{equation}
and $\tilde \varphi'$ is a twist of $\tilde \varphi$ by an unramified character of $F^\times$ such that $\tilde \varphi'$ is also self-dual and $\tilde \varphi'\not\cong\tilde \varphi$. If $G$ is unramified unitary, then we just take \begin{equation}
\label{unitary-parameter}
\varphi = \bigoplus_{\tilde\varphi} \tilde\varphi[[a_{\tilde\varphi}]]\oplus \tilde \varphi'[[b_{\tilde\varphi}]],
\end{equation}
with $\tilde\varphi $ in the sum ranging over the same collections (\ref{self-dual-condition-deg>1}), but without requiring $[E:F]>1$. There are extra conditions on the numbers $a,b,c,d,(a_{\tilde\varphi},b_{\tilde\varphi})_{{\tilde\varphi}}$ for $\varphi$ to have image in ${}^L\mathbf G$.
\begin{enumerate}[(i)]


\item If $G$ is symplectic or even orthogonal, then $a,b,c,d$ are all odd with extra conditions. 

In the former case, since $\deg\varphi$ is odd, three of these $a,b,c,d$ are $\equiv i \bmod 4$ where $i=1$ or 3, and the remaining one is $\equiv 4-i \bmod 4$. Since $\det\varphi =\omega_0^{(b+1)/2}\omega_1^{(c+1)/2}\omega_2^{(d+1)/2} =\omega_0^{(b+d)/2+1} \omega_1^{(c+d)/2+1} =  1$, the three are $b,c,d$ that have the same congruence mod 4. 

In the latter case, since $\deg\varphi$ is even, either all $a,b,c,d$ are $\equiv i \bmod 4$ where $i=1$ or 3, in which cases $\det \varphi =1$, or exactly two of them are either $\equiv 1$ or  $\equiv 3 \bmod 4$, in which cases $\det \varphi =\omega_i$ if $1[[a]]\oplus \omega_i[[a_i]]\hookrightarrow \varphi$ with $a_i\equiv a\mod 4$. The character $\det \varphi$ determines the quadratic field extension $F[\sqrt d]/F$, where $d$ is the discriminant of the defining Hermitian form $h_V$ of $G=G(V)$.

\item For all other $(a_{\tilde\varphi},b_{\tilde\varphi})_{{\tilde\varphi}}$, since $\tilde\varphi$ and $\tilde\varphi'$ are of opposite parities, so are the pair $a_{\tilde\varphi}$ and $b_{\tilde\varphi}$. If $G$ is orthogonal or symplectic, then $a_{\tilde\varphi}$ is odd; if $G$ is unitary, then $a_{\tilde\varphi}$ is odd if and only if $N+\deg \tilde\varphi$ is even.

\end{enumerate}
We remark that the exceptional case for even orthogonal groups happens when  $a,b,c,d$ are all $-1$.

We sometimes shorthand denote a parameter $\varphi $ for $G$ as 
\begin{equation}
\label{short-form-parameter}
\varphi = [[a,b,c,d,(a_{\tilde\varphi},b_{\tilde\varphi})_{{\tilde\varphi}}]]\text{ for $\varphi $ in (\ref{quadratic-parameter}),\quad or }\quad
[[(a_{\tilde\varphi},b_{\tilde\varphi})_{{\tilde\varphi}}]]\text{ for $\varphi $ in (\ref{unitary-parameter})}.
\end{equation}
Take another parameter $\varphi' = [[a',b',c',d',(a'_{\tilde\varphi},b'_{\tilde\varphi})_{{\tilde\varphi}}]]$, where all $a',b',c',d',(a'_{\tilde\varphi},b'_{\tilde\varphi})$ satisfy the above conditions. We call $\varphi'$ a {\bf companion} of $\varphi$ if $\{a,b\}= \{a',b'\}$, $\{c,d\}= \{c',d'\}$, and $(a_{\tilde\varphi},b_{\tilde\varphi}) = (a'_{\tilde\varphi},b'_{\tilde\varphi})$ for all $\tilde\varphi$. Hence a parameter has at most 4 companions (including itself).

Define the component group $A(\varphi):=\pi_0(Z_{{\hat{\mathbf G}}}(image\varphi))$. Given a parameter $\varphi$ as in (\ref{quadratic-parameter}) or (\ref{unitary-parameter}), denote by ${r(\varphi)}$ the number of irreducible components of $\varphi$, then $A(\varphi)$ is isomorphic to 
\begin{equation*}
\begin{split}
&\{\pm 1\}^{r(\varphi)-1}\text{ when $G$ is symplectic;}
\\
&\{\pm 1\}^{r(\varphi)}\text{ when $G$ is even orthogonal or unitary.}
\end{split}
\end{equation*}
We also define $A^0(\varphi)$ be the image of $Z_{{\hat{\mathbf G}^0}}(image\varphi)$ in $A(\varphi)$ when $G$ is even orthogonal, in which case $A^0(\varphi)\cong \{\pm 1\}^{r(\varphi)-1}$.

Since we will not go into the study of the genericity of representations, while it is not too difficult to give the definition in \cite{Moeglin-multiplicity-1, Moeg-base-change} of alternating characters of $A(\varphi)$, we only provide the number of such characters for the given parameter $\varphi$ as in (\ref{quadratic-parameter}) or (\ref{unitary-parameter}). Let $k(\varphi)$ be the number of  elements in the following multi-sets:
\begin{equation*}
\begin{split}
\{a,b,c,d, \{a_{\tilde\varphi},b_{\tilde\varphi}\}_{\tilde\varphi}\}\cap (2\mathbb N-1) \quad &\text{when $G$ is orthogonal or symplectic};
\\
\{a_{\tilde\varphi},b_{\tilde\varphi}\}_{\tilde\varphi}\cap (2\mathbb N-1) \quad &\text{when $G$ is unitary.}
\end{split}
\end{equation*}
The number of alternating characters of $\varphi$ is equal to
\begin{equation}
\label{number of alternating characters}
\begin{split}
&2^{k(\varphi)-1}\text{ when $G$ is symplectic};
\\
&2^{k(\varphi)}\text{ when $G$ is even orthogonal or unitary.}
\end{split}
\end{equation}
In the rest of the paper, when we say a parameter for $G$, we mean a parameter of the form in (\ref{quadratic-parameter}) or (\ref{unitary-parameter}).

\subsection{LLC for classical groups}
\label{section LLC for classical groups}

We first recall some facts about the LLC for depth zero supercuspidal representations of general linear groups, with an explicit description. The LLC for $\tilde G = \GL(m)$ asserts that supercuspidal representations of $\tilde {G}$ can be parametrized by irreducible representations of $\mathcal W_F$ of degree $m$. If a parameter $\tilde\varphi$ is trivial on the wild-inertia subgroup of $\mathcal W_F$, so that $\tilde\varphi \cong \Ind^{\mathcal W_F}_{\mathcal W_E}\tilde\xi$, where $E/F$ is the unramified extension of degree $m:=\deg\tilde\varphi$, and $\tilde\xi$ is a regular character of $E^\times$, then the representation $\tilde\pi = \tilde\pi_{\tilde\varphi}$ corresponding $\tilde\varphi$ under the LLC can be constructed directly from $\tilde\xi$ as follows.

We first descend $\tilde\xi|_{\mathcal O_E^\times}$ to a character $\overline{\tilde\xi}$ of the group $\GL(1,\mathbb F_{q^m})$ viewed as an elliptic torus in $ {\tilde{\mathsf G}} =\mathsf{GL}(m)$. If $\overline{\tilde\xi}$ corresponds to an element  $\tilde s\in \hat{\mathsf T} \subseteq \hat {\tilde{\mathsf G}}$, then its regularity implies that $\tilde s$ corresponds to a cuspidal representation $\tilde\rho_{\overline{\tilde\xi}}\in \mathcal E(\tilde s)$. In fact, $\tilde\rho_{\overline{\tilde\xi}}$ is the Deligne-Lusztig induction of ${\overline{\tilde\xi}}$. We then inflate $\tilde\rho_{\overline{\tilde\xi}}$ to $\GL(m,\mathcal O_F)$. By \cite[2.4 Th 2.]{BH-ET3}, if $\tilde{\xi}(\varpi)=(-1)^{m-1}$, and if $\tilde{\rho}_{\tilde{\xi}}:={{\tilde\xi}}\tilde\rho_{\overline{\tilde\xi}}$ is the inducing cuspidal type on $\tilde{\mathcal J} $ as constructed in Section \ref{subsection Compact inductions}, then the supercuspidal representation $
\tilde \pi=\cInd_{ \tilde{\mathcal J}}^{\tilde G}\tilde{\rho}_{\tilde{\xi}}$ as constructed in (\ref{depth zero supercuspidal compact induction}) corresponds to $\tilde\varphi$ under the LLC.

We now turn to gather together some facts about LLC for classical groups. Suppose first that $G=\mathbf G(F)$ is a classical group. Let $\mathbf G_-$ be the `pure inner form' of $\mathbf G_+=\mathbf G$, i.e., $\mathbf{G}_+$ and $\mathbf{G}_-$ are parametrized by the cohomology group $H^1(\Fo,\mathbf G)$, which is either trivial (in which case $\mathbf G_+=\mathbf G_-$) or isomorphic to $\{\pm 1\}$ for classical groups. The LLC for a classical group $G^0$ is a canonical surjection
$$\mathrm{LLC}:\bigsqcup_{\delta\in \{\pm \}}\left\{\begin{matrix}
\text{isomorphism classes}
\\
\text{of irreducible}
\\
\text{representations of $G^0_\delta$}
\end{matrix}\right\} \rightarrow 
\{
\text{Langlands parameters of $G^0$}
\}$$
satisfying a number of canonical properties, among those the stability: if $\Pi_\varphi$ is the preimage under the LLC of $\varphi$, then members in $\Pi_\varphi$ satisfy the twisted endoscopic identity. 

For odd orthogonal, symplectic, and unitary groups, the existence of the LLC is proved by \cite{Arthur-new-book, Mok-unitary, unitary-inner}. For even orthogonal groups: while the LLC for disconnected groups is still underdeveloped, it is actually `more natural', from Arthur's theory \cite{Arthur-new-book}, to establish the LLC for the full even orthogonal group $G$ rather than $G^0$. This is in fact the view-point of Arthur: simply speaking, he established a canonical bijection 
$$\mathrm{LLC}:\bigsqcup_{\delta\in \{\pm \}}\left\{\begin{matrix}
\text{isomorphism classes}
\\
\text{of irreducible}
\\
\text{representations of $G^0_\delta$}
\end{matrix}\right\}/G_\delta\text{-conjugacy} \rightarrow 
\{
\text{Langlands parameters of $G^0$}
\}
/\hat{\mathbf G}\text{-conjugacy}.$$
Note that parameters on the RHS above can be viewed as isomorphism classes of orthogonal representations of the Weil-Deligne group of an even degree. We certainly have the natural restriction map on the representation side
$$
\left\{\begin{matrix}
\text{isomorphism classes}
\\
\text{of irreducible}
\\
\text{representations of $G$}
\end{matrix}\right\}/G\text{-conjugacy} \rightarrow 
\left\{\begin{matrix}
\text{isomorphism classes}
\\
\text{of irreducible}
\\
\text{representations of $G^0$}
\end{matrix}\right\}/G\text{-conjugacy}.$$
However, it is in general unknown whether the surjective map from $\hat{\mathbf G}$- to $\hat{\mathbf G}^0$-conjugacy classes of parameters is naturally compatible with this restriction. It could be possible to solve this problem by studying the $G/G^0$-torsor structure on both sides of the LLC as in \cite{Arthur-even-orthogonal}.

The representations in $\Pi_{\varphi}=\Pi_{\varphi}(G_+)\cup \Pi_{\varphi}(G_-)$ can then be parametrized by the Pontryagin dual $A(\varphi)^\wedge$ of $A(\varphi)$. We assume that a generic representation (i.e., possessing a  Whittaker model) in $\Pi_{\varphi}$ is chosen to be parametrized by the trivial character of $A(\varphi)$. For $\epsilon\in \{\pm\}$, we have $\pi \in\Pi_{\varphi}(G_\epsilon)$ if and only if
\begin{equation}
\label{Nakayama-Tate}
\text{the restriction of its corresponding character to $Z\hat{\mathbf G}$ is $\epsilon\in (Z\hat{\mathbf G})^\wedge \cong H^1(F_\bullet,\mathbf G)$}.
\end{equation}
Moreover, $\Pi_{\varphi}$ contains a, and in fact consists of, discrete series if and only if it is a direct sum of multiplicity-free irreducible components, c.f., \cite{Moeg-exhaustion, Moeg-base-change, Moeglin-Renard-non-quasi-split}. Our parameters of the form $\varphi = [[a,b,c,d,(a_{\tilde\varphi},b_{\tilde\varphi})_{{\tilde\varphi}}]]$ in (\ref{short-form-parameter}) satisfy this condition.

According to \cite[Th 1.5.1]{Moeglin-multiplicity-1}, \cite[Th 8.4.4]{Moeg-base-change}, the supercuspidal representations in $\Pi_\varphi = \Pi_\varphi(G_+) \cup \Pi_\varphi(G_-) $ are parametrized by the alternating characters in $A(\varphi)^\wedge$. The numbers are given by (\ref{number of alternating characters}), and let's reiterate: the number of supercuspidals in $\Pi_\varphi$ is
\begin{equation}
\label{number of supercuspidals in a packet}
\begin{split}
&2^{k(\varphi)-1}\text{ when $G$ is symplectic};
\\
&2^{k(\varphi)}\text{ when $G$ is even orthogonal or unitary}.
\end{split}
\end{equation}
For even special orthogonal groups, the number of supercuspidals in $\Pi_\varphi$ for a parameter $\varphi$ of $G^0$ is counted using $A^0(\varphi)$ instead of $A(\varphi)$. In the exceptional case when $\varphi$ has all irreducible components of even dimension, we have to count each representation twice. Hence the number of supercuspidals in $\Pi_\varphi$ is
\begin{equation}
\label{number of supercuspidals in a packet even SO}
\begin{split}
&2^{k(\varphi)-1}\text{ in the non-exceptional case};
\\
&2^{k(\varphi)}\text{ in the exceptional case.}
\end{split}
\end{equation}
Combining with (\ref{Nakayama-Tate}), we can precisely count the number of supercuspidal representations in an individual $\Pi_\varphi(G_\epsilon)$.

Finally, we provide a result of M{\oe}glin linking the information on the irreducible components of parameters with the reducibility of parabolic inductions for classical groups. Let $\tilde\pi = \tilde\pi_{\tilde\varphi}$ be a supercuspidal representation of $\tilde G=\GL(m)$ with parameter $\tilde\varphi$, and $\pi\in \Pi_{\varphi}$ is a supercuspidal representation of $G^0$ lying in the packet corresponding to $\varphi$. If $\tilde\varphi[[a]]$ is a component of a parameter $\varphi$, the number $a=a_{\tilde\varphi}$ is closely related to the points where the following parabolically induced representation is reducible:
\begin{equation}
\label{parabolically induced representation}
I(s,\tilde\pi,\pi):= \iota_{P}^{G^0_W}(|\det|^s\tilde\pi_{}\times \pi),\qquad s\in \mathbb R
\end{equation}
(for any parabolic subgroup $P$ in $G^0_W$ containing a Levi subgroup of the form $\tilde G\times G^0$). By \cite{silberger-special}, there are at most two such points of the form $\pm s$ for some $s\in \mathbb R_{\geq 0}$. By M{\oe}glin's theory in \cite[Sec 1.3, 1.4]{Moeglin-multiplicity-1}, \cite[Sec 3, 6]{Moeg-base-change}, we have
\begin{equation}
\label{SL2-component-and-reducibility}
\text{
$\tilde\varphi[[a]]$ is a component of $\varphi$ \quad if and only if  \quad $I(s,\tilde\pi,\pi)$ is reducible at $s=(a+1)/2$.}
\end{equation}
We will more generally consider points of reducibility of $I(s,\tilde\pi,\pi)$ for $s\in \mathbb C$. For a fixed representation $\tilde\pi$
as in (\ref{depth zero supercuspidal compact induction}), the group of unramified characters  of $\tilde G = \GL(m)$ stabilizing $\tilde\pi$ by tensoring is cyclic of order $m$. When we talk about points of reducibility later on, we assume that the points are lying in the domain \begin{equation}
\label{fundamental domain of reducibility}
s\in \mathbb C,\quad \text{ where}\quad 0\leq Im(s)< \frac{2\pi }{m\log q},
\end{equation}
which determines the whole (countably infinite) set of points of reducibility in $\mathbb C$.

\subsection{Main results on typically almost symmetric representations}
\label{subsection The main result}

We now provide our condition on the depth zero supercuspidal representations we want to study. Let $\rho = \rho (m_+,m_-,(m_{\tilde P})_{\tilde P})$ be a cuspidal representation of a classical group. We call another such representation $\rho' = \rho (m'_+,m'_-,(m'_{\tilde P})_{\tilde P})$ a {\bf companion} of $\rho$ if 
\begin{equation}
\label{condition_rhoy_equals_rhoz}
\text{$|m_{\tilde P}| = |m'_{\tilde P}|$ \quad  for all $\tilde P$.}
\end{equation}
If $\pi$ is a supercuspidal representation of a classical group $G$ and is compactly induced from a type constructed from $\rho_y\times \rho_z$, we say that $\pi$ is {\bf typically almost symmetric} if $\rho_y$ is a companion of $\rho_z$, and is {\bf typically symmetric} if $\rho_y = \rho_z$. This condition restricts our attention to $G$ being symplectic, even orthogonal, and even unramified unitary, as $\mathsf G_y = \mathsf G_z$.

Later after Section \ref{subsection Calculations with Hecke algebras}, when we develop respectively a relation between the pairs of tuples 
\begin{equation*}
\begin{split}
&(m_{+,y},  m_{+,z})
\quad \text{ and }\quad (a,b);
\\
&(m_{-,y},  m_{-,z})
\quad \text{ and }\quad (c,d);
\\
&(m_{\tilde P,y},  m_{\tilde P,z})
\quad \text{ and }\quad (a_{\tilde \varphi}, b_{\tilde \varphi}),\text{ for all other $\tilde \varphi$}, 
\\
&\qquad\qquad \text{where $\tilde{\xi}|_{\mathcal O_E^\times}$ corresponds to a root of $\tilde{P}$ and $\tilde{\xi}(\varpi)=1$,}
\end{split}
\end{equation*}
via parameters of Hecke algebras and M{\oe}glin's relation (\ref{SL2-component-and-reducibility}), then we will see that the typically almost symmetry on $\pi$ is equivalent to that 
one of the numbers in each tuple $(a,b)$ and $(c,d)$, and the odd number in each $(a_{\tilde \varphi}, b_{\tilde \varphi})$, is $-1$. As a more precise result, we state our first main proposition on reducibility.

\begin{prop}
\label{main results on reducibility points}
Fix a depth zero supercuspidal representation $\pi$ which is moreover typically symmetric, i.e., $\pi$ is constructed from a type $\rho_y\times \rho_z$ with $\rho_y=\rho_z$. Suppose both $\rho_w =\rho(m_+,m_-,(m_{\tilde P})_{\tilde P})$ or $\rho((m_{\tilde P})_{\tilde P})$. We list the set $\Red(\tilde\pi,\pi)$ of points of reducibility for $I(s,\tilde\pi,\pi)$ (\ref{parabolically induced representation}), with $\tilde\pi$ constructed from a self-dual polynomial $\tilde P$.
\begin{enumerate}[(i)]

\item When either $G$ is even unramified unitary, or $\tilde\pi\not\in\{1,\omega_0,\omega_1,\omega_2\}$ when $G$ is symplectic or even orthogonal. If  $\tilde\pi = \tilde\pi_{\tilde\rho}$ with $\tilde\rho(\varpi)=1$, then 
$$\Red(\tilde\pi,\pi)  = \left\{ \pm \left( m_{\tilde P} + \frac{1}{2}\right), \pm 0+\frac{\pi i}{\deg \tilde P\log q}\right\}.$$
\label{main results on reducibility points even unitary}

\item When $G$ is 
symplectic, then 
$$\Red(\omega_1,\pi)  = \left\{ \pm 2m_{-}, \pm 0+{\pi i}/{\log q}\right\}.$$
The same statement holds when $G$ is even orthogonal and both $\mathsf G_w$ are even, in which case we have furthermore
$$\Red(1,\pi)  = \left\{ \pm 2m_{+}, \pm 0+{\pi i}/{\log q}\right\}.$$
\label{main results on reducibility points symplectic}

\item When $G$ is even orthogonal and both $\mathsf G_w$ are odd, then
$$\Red(1,\pi)  = \left\{ \pm (2m_{+}+1), \pm 0+{\pi i}/{\log q}\right\} \quad\text{and}\quad
\Red(\omega_1,\pi)  = \left\{ \pm (2m_{-}+1), \pm 0+{\pi i}/{\log q}\right\}.$$

\end{enumerate}
When $G$ is even orthogonal, both $\mathsf G_w$ are even, and $\rho_y = \Ad(h)\circ \rho_z$, where $h\in \mathsf{GG(V)}$ is defined before Lemma \ref{properties as in Waldspurger book}, then 
$$\Red(1,\pi)  = \left\{ \pm 2m_{+}, \pm 0+{\pi i}/{\log q}\right\} \quad\text{and}\quad
\Red(\omega_1,\pi)  = \left\{ \pm 0 , \pm 2m_{-}+{\pi i}/{\log q}\right\}.$$
The $\pm 0$ notations emphasize that we always have four points of reducibility in the domain (\ref{fundamental domain of reducibility}), some of which coincide with each other when Re$(s)=0$ for some $s\in \Red(\tilde\pi,\pi)$. 
\end{prop}
Note that we only list the reducibility results that are essential to our final arguments in Section \ref{subsection The final argument}, as we will assume (\ref{number of supercuspidals in a packet}) on the expected number of supercuspidal representations in a packet to avoid repeating calculations and facilitate the arguments.

Combining Proposition \ref{main results on reducibility points} with (\ref{SL2-component-and-reducibility}), we can determine the components of the parameter $\varphi$ whose packet contains $\pi$. We first state our first result of this sort for unitary groups.

 \begin{prop}
  \label{main result for unramified unitary}
Suppose that $\tilde\varphi[[a_{\tilde\varphi}]]\oplus \tilde \varphi'[[b_{\tilde\varphi}]]$ appears as a component of $\varphi$ in the cases of either

\begin{enumerate}[(i)]

\item \unskip $G$ is even unramified unitary, or 

\item \unskip $\tilde\varphi\not\in\{1,\omega_0,\omega_1,\omega_2\}$ and $G$ is symplectic or even orthogonal.

\end{enumerate}
\setlength{\parskip}{0pt}
Then $a_{\tilde\varphi} = 2m_{\tilde P}$ and $b_{\tilde\varphi} = -1$. 
\end{prop}

As a consequence, we can determine the parity of a depth zero cuspidal representation of a general linear group.

\begin{cor}
\label{Corollary for parity on representation side}
The representation $\tilde\pi = \cInd_{\tilde{\mathcal J}}^{\tilde G}\tilde \rho$ of $\tilde G = \GL(m)$ constructed in (\ref{depth zero supercuspidal compact induction}), with $\tilde{\rho}|
_{\mathcal O_E^\times }$ conjugate-self-dual and $\tilde{\rho}(\varpi)=1$, is conjugate-orthogonal (resp. conjugate-symplectic) if  $m$ is odd (resp. even).
\end{cor}

One may obtain the same result in Corollary \ref{Corollary for parity on representation side} via the explicit LLC for general linear groups, see Section \ref{(Comparison with the LLC for GL(n))}. Our proof in this article is almost completely on the representation side.

We then state another result which can be derived from Proposition \ref{main results on reducibility points}: the classification of typically almost symmetric depth zero supercuspidal representations into L-packets. We remark that the result refines an ambiguity issue in \cite{Lust-Stevens}, where a family of depth zero supercuspidal representations, including the typically almost symmetric ones, can only be classified into a finite union of multiple packets.

Take a typically almost symmetric representation $\pi = \pi_\rho$ of $G$ constructed from $\rho_y\times \rho_z$, with non-negative integers $|m_+|$,  $|m_-|$, and all $(m_{\tilde P})_{\tilde P\neq X\pm 1}$ fixed.

 \begin{prop}
\label{main-result-symplectic}
Suppose that $G$ is a symplectic group, and put $a=4m_++1$ and $c=4m_--1$, then  
\begin{equation*}
\text{
$\pi_\rho\in \Pi_{[[a,-1,c,-1,(2m_{\tilde P}, -1)]]}$\quad  if and only if\quad  $\rho_y = \rho_z = \rho(m^+, \pm m^-,(m_{\tilde P})_{\tilde P})$.
}
\end{equation*}
and 
\begin{equation*}
\begin{split}
\pi_\rho\in \Pi_{[[a,-1,-1,c,(2m_{\tilde P}, -1)]]}\quad &\text{ if and only if}\quad  
\\
&\{\rho_y ,\rho_z\} = \{\rho(m^+,  m^-,(m_{\tilde P})_{\tilde P}),\rho(m^+, -m^-,(m_{\tilde P})_{\tilde P})\}.
\end{split}
\end{equation*}
\end{prop}
Note that, we always have $b\equiv c\equiv d \bmod 4$ for parameters of symplectic groups, and typically almost symmetry implies that one of $c$ and $d$ is $-1$, hence we must have $b=-1$. Moreover, under the assumption stated in Proposition \ref{main-result-symplectic}, $c=-1$ if and only if $m_-=0$, and in which case we only have one supercuspidal induced by $\rho_y=\rho_z = \rho(m^+, 0,m_{\tilde P})$, the only supercuspidal member in $\Pi_{[[a,-1,-1,-1,(2m_{\tilde P}, -1)]]}$.

The results for even orthogonal groups is similar, with a bit modification.
 \begin{prop}
 \label{main-result-orthogonal}
Suppose that $G$ is an even orthogonal group and both $\mathsf G_w$, where $w\in \{y,z\}$, are also even orthogonal. Put $a = 4m_+ -1 $ and $c = 4m_--1 $, then 
\begin{equation*}
\begin{split}
\pi_\rho\in \Pi_{[[a,-1,c,-1,(2m_{\tilde P}, -1)]]} \quad &\text{ if and only if }\quad  \rho_y=\rho_z = \rho(\epsilon m^+, \epsilon' m^-, (m_{\tilde P})_{\tilde P}),
\\
\pi_\rho\in \Pi_{[[-1,a,-1,c,(2m_{\tilde P}, -1)]]} \quad &\text{ if and only if }\quad  
\\
&\{\rho_y,\rho_z\}= \{\rho(\epsilon m^+, \epsilon' m^-,(m_{\tilde P})_{\tilde P}),\rho(-\epsilon m^+, -\epsilon ' m^-,(m_{\tilde P})_{\tilde P})\},
\\
\pi_\rho\in \Pi_{[[a,-1,-1,c,(2m_{\tilde P}, -1)]]} \quad &\text{ if and only if }\quad  
\\
&\{\rho_y,\rho_z\}= \{\rho(\epsilon  m^+,\epsilon' m^-,(m_{\tilde P})_{\tilde P}),\rho(\epsilon  m^+, -\epsilon' m^-,(m_{\tilde P})_{\tilde P})\},
\\
\pi_\rho\in \Pi_{[[-1,a,c,-1,(2m_{\tilde P}, -1)]]} \quad &\text{ if and only if }\quad  
\\
&\{\rho_y,\rho_z\}= \{\rho(\epsilon m^+, \epsilon' m^-,(m_{\tilde P})_{\tilde P}),\rho(-\epsilon m^+, \epsilon' m^-,(m_{\tilde P})_{\tilde P})\},
\end{split}
\end{equation*}
with $\epsilon,\epsilon'\in\{\pm \}$.
\end{prop}
We leave the analyses of the situations for $a=-1$ or $c=-1$ as exercises.

 \begin{prop}
 \label{main-result-orthogonal-odd-odd}
Suppose that $\mathsf G_w$, where $w\in \{y,z\}$, are odd orthogonal. Put $a = 4m_+ +1 $ and $c = 4m_-+1 $, then for $G$ the even orthogonal group with discriminant $d=\varpi$, 
\begin{equation*}
\begin{split}
\pi_\rho\in \Pi_{[[a,-1,c,-1,(2m_{\tilde P}, -1)]]} \quad &\text{ if and only if }\quad  \rho_y=\rho_z = \rho( m^+, m^-, (m_{\tilde P})_{\tilde P})(\epsilon),
\\
\pi_\rho\in \Pi_{[[-1,a,-1,c,(2m_{\tilde P}, -1)]]} \quad &\text{ if and only if }\quad  
\\
&\{\rho_y,\rho_z\}= \{\rho( m^+, m^-, (m_{\tilde P})_{\tilde P})(\epsilon),\rho( m^+, m^-, (m_{\tilde P})_{\tilde P})(-\epsilon)\},
\end{split}
\end{equation*}
with $\epsilon \in \{ \pm \}$ and for each pure inner form of $G$. 

The same result holds for $G$ the even orthogonal group with discriminant $d=\zeta\varpi$, where $\zeta\in \mu_F \smallsetminus \mu_F^2$, with the packets $\Pi_{[[a,-1,c,-1,(2m_{\tilde P}, -1)]]}$ and $\Pi_{[[-1,a,-1,c,(2m_{\tilde P}, -1)]]} $ respectively replaced by  $\Pi_{[[a,-1,-1,c,(2m_{\tilde P}, -1)]]}$ and $\Pi_{[[-1,a,c,-1,(2m_{\tilde P}, -1)]]} $. \end{prop}

Finally, we provide the classification for even special orthogonal groups as a record.

\begin{prop} 
\label{main-result-connectd-orthogonal}
Continue from Proposition \ref{main-result-orthogonal} or \ref{main-result-orthogonal-odd-odd}, let $G^0$ be the subgroup of the connected component of $G$. 
\begin{enumerate}[(i)]
\item
When both $\mathsf G_w$ are odd, or when both $\mathsf G_w$ are even and $\varphi$ is non-exceptional, i.e., $(a,c)\neq (-1,-1)$, let $\varphi^0$ be the parameter of $G^0$ such that $\varphi = \varphi^0$ mod $\hat{\mathbf G}$-conjugacy. Let $\rho^0 = \rho(\pm)$ as defined in (\ref{(++)on-S(OtimesO)-odd-odd}) or (\ref{(++)on-S(OtimesO)}), then 
$$\pi_{\rho^0}\in \Pi_{\varphi^0}\quad\text{if and only if}\quad \pi_{\rho}\in \Pi_{\varphi}.$$

\item
When both $\mathsf G_w$ are even and $\varphi$ is exceptional, i.e.,  $a=c=-1$, for $\epsilon\in\{\pm\}$, let $\varphi_\epsilon$ be the two parameters of $G^0$ such that $\varphi = \varphi_\epsilon$ mod $\hat{\mathbf G}$-conjugacy, and let $\rho^0 = \rho_\epsilon$, as defined in (\ref{all-zero-on-S(OtimesO)}), then 
$$\Pi_{\varphi_+}\sqcup \Pi_{\varphi_-}\quad\text{contains exactly two supercuspidal representations}\quad \pi_{\rho_+}\text{ and }\pi_{\rho_-}.$$
\end{enumerate}
\end{prop}
Note that the union of two packets in the exceptional case is the best we can do, according to the theory of endoscopy for even orthogonal groups, c.f. \cite{Arthur-new-book, Arthur-even-orthogonal, atobe-Gan-SO-even}.

The proof of all these propositions are similar. They will be given in Section \ref{subsection The final argument}, after we acquire certain knowledge in Hecke algebras for covering types in the next chapter.

\section{Hecke algebras of covering types}

After providing the preliminary knowledge on the relation between Bernstein components and modules over Hecke algebras across types and their covers  for general reductive groups (Sec \ref{section Covering Types}), we recall the construction of covering types for classical groups (Sec \ref{subsection Covering types for classical groups}) and describe their associated Hecke algebras and modules (Sec \ref{subsection Hecke algebras and modules}). Viewing elements in the Hecke algebras as operator-valued maps, we introduce a notion of canonical normalizations on the intertwining operators associated to the generators of the Hecke algebras (Sec \ref{subsection Canonical normalization of intertwining operators}), which are used to study the reducibility of modules over Hecke algebras via the key relation (\ref{Final comparison of eigenvalues}) deduced from Proposition \ref{Blondel-Blasco about product of eigenvalues} on comparing eigenvalues. Further calculations in Section \ref{subsection Proof of the main result} are required to prepare our arguments in Chapter \ref{subsection The final argument} for the proofs of our main results. Section \ref{subsection Canonical normalization of intertwining operators} serves as an appendix to prove a non-vanishing result of the character of the intertwining operators.

\subsection{Covering types}
\label{section Covering Types}

 Let $\mathcal{M}$ be a connected reductive group, and $\pi_{\mathcal{M}}$ be an irreducible representation of $\mathcal{M}$. Denote by $\mathfrak{s}_{\mathcal{M}}$ the inertial class of $\pi_{\mathcal{M}}$. Let $\lambda_{\mathcal{M}}$ be an irreducible representation of a compact open subgroup 
 $\mathcal{J}_{\mathcal{M}}$ of $\mathcal{M}$. We call $(\mathcal{J}_{\mathcal{M}},\lambda_{\mathcal{M}})$ an $
 \mathfrak{s}_{\mathcal{M}}$-type if the functor 
 \begin{equation}
 \label{First equivalence on types}
\bf M_{\mathcal M}:\mathcal{R}^{\mathfrak{s}_{\mathcal{M}}}({\mathcal{M}})\rightarrow \text{Mod-}\mathcal{H}({\mathcal{M}},\lambda_{\mathcal{M}} ),\, \qquad V\mapsto \Hom_{\mathcal{J}_{\mathcal{M}}}(V_{\lambda_{\mathcal{M}} },V),
 \end{equation}
between
 \begin{itemize}
 \item $\mathcal{R}^{\mathfrak{s}_{\mathcal{M}}}({\mathcal{M}}) $ the Bernstein component of $\mathfrak{s}_{\mathcal{M}}$, i.e., the sub-category of representations 
 of $\mathcal{M}$, all of whose objects having irreducible 
sub-quotients with cuspidal support in $\mathfrak{s}_{\mathcal{M}}$, and

 \item $\text{Mod-}\mathcal{H}({\mathcal{M}},\lambda_{\mathcal{M}})$ the category of non-degenerate right modules over the ${\mathcal{M}}$-intertwining algebra $\mathcal{H}({\mathcal{M}},\lambda_{\mathcal{M}})$ of $\lambda_{\mathcal{M}}$,
 \end{itemize}
 is an equivalence of categories. In particular, suppose that $\pi_\mathcal{M}$ is a supercuspidal representation of the form 
$\cInd_{\mathbf{J}_\mathcal{M}}^{\mathcal{M}}\boldsymbol{\lambda}_{\mathcal{M}}
$,  where $\boldsymbol{\lambda}_{\mathcal{M}}$ is an irreducible representation of an open compact-mod-center subgroup of $\mathcal{M}$. Let $\mathcal{J}_{\mathcal{M}}$ be the unique  maximal compact subgroup of $\mathbf{J}_{\mathcal{M}}$ and $\lambda_{\mathcal{M}}$ be an irreducible component of $\boldsymbol{\lambda}_{\mathcal{M}}|_{\mathcal{J}_{\mathcal{M}}}$, then \cite[(5.4)]{BK-cover} $(\mathcal{J}_{\mathcal{M}},\lambda_{\mathcal{M}})$ is an $
 \mathfrak{s}_{\mathcal{M}}$-type. Note that there is an obvious action of the group of unramified characters of $\mathcal M$ on both sides of (\ref{First equivalence on types}) such that $\bf{M}_{\mathcal{M}}$ is equivariant.

Suppose now that ${\mathcal{M}}$ is a Levi subgroup of another connected reductive group ${\mathcal{G}}$, and $\mathcal J$ is a compact open subgroup of ${\mathcal{G}}$. An irreducible representation $(\mathcal{J},\lambda)$ is called a {\bf cover} of $(\mathcal{J}_{\mathcal{M}},\lambda_{\mathcal{M}})$ if, for a \underline{given} parabolic (see \cite[Comments before (8.8)]{BK-cover}) 
 ${\mathcal{P}}={\mathcal{M}}{\mathcal{N}}$ with Levi component ${\mathcal{M}}$, unipotent radical ${\mathcal{N}}= {\mathcal{N}}_+$, and the opposite ${\mathcal{N}}_-$ of $\mathcal N$,
 \begin{enumerate}[(i)]
 \item $\mathcal{J}$ decomposes with respect to ${\mathcal{M}}$, with $\mathcal{J}\cap {\mathcal{M}} = \mathcal{J}_{\mathcal{M}}$;
 \item $\lambda|_{\mathcal{J}_{\mathcal{M}}} = \lambda_{\mathcal{M}}$, and both $\lambda|_{\mathcal{J}_+}$ and $\lambda|_{\mathcal{J}_-}$ lie in the kernel of $\lambda$, where $\mathcal{J}_{\pm }:= \mathcal{J}\cap {\mathcal{N_{\pm }}}$;
 
 \item there exists an invertible element in  $\mathcal{H}({\mathcal{G}},\lambda)$ supported on a $\mathcal{J}$-double coset $\mathcal{J}\varpi_{\mathcal{P}}\mathcal{J}$, where $\varpi_{\mathcal{P}}$ is strongly positive with respect to $({\mathcal{P}},\mathcal{J})$, in the sense of \cite[(6.16)]{BK-cover}. 
 \end{enumerate}
 In particular, we have the Iwahori decomposition for $\mathcal{J}$ with respect to ${\mathcal{P}}={\mathcal{M}}{\mathcal{N}}$, i.e., $\mathcal{J} = \mathcal{J}_+\mathcal{J}_{\mathcal{M}}\mathcal{J}_-$. Then by \cite[(8.3)]{BK-cover}, under the above conditions, $(\mathcal{J}_{\mathcal{P}},\lambda_{\mathcal{P}})$ is an $\mathfrak{s}$-type where $\mathfrak{s}:=[\mathfrak{s}_{\mathcal{M}}]_{\mathcal{G}}$ is the inertial class ${\mathcal{G}}$ with support $\mathfrak{s}_{\mathcal{M}}$. This leads to the following commutative diagram
\begin{equation}\label{commutative diagram}
  \xymatrixcolsep{5pc}\xymatrix{
 \mathcal{R}^{\mathfrak{s}}({\mathcal{G}})\ar[r]^{\bf M_{\mathcal G}}    
&\text{Mod-}\mathcal{H}({\mathcal{G}},{\lambda_{\mathcal{P}}})
\\
 \mathcal{R}^{\mathfrak{s}_{\mathcal{M}}}({\mathcal{M}})\ar[u]_{\iota_P^{{\mathcal{G}}}}
\ar[r]^{\bf M_{\mathcal M}}
&\text{Mod-}\mathcal{H}({\mathcal{M}},\lambda_{\mathcal{M}}),  
\ar[u]_{(t_{\mathcal{P}})_*}
}
\end{equation}
where $\iota_{\mathcal{P}}^{\mathcal{G}}$ is the normalized parabolic induction functor, and $(t_{\mathcal{P}})_*$ is the functor naturally induced from the embedding 
\begin{equation*}
t_{\mathcal P}:\mathcal{H}({\mathcal{M}},\lambda_{\mathcal{M}}) \rightarrow \mathcal{H}({\mathcal{G}},\lambda_{\mathcal{P}})
\end{equation*}
defined by \cite[(7.12)]{BK-cover}. 

\subsection{Covering types for classical groups}
\label{subsection Covering types for classical groups}


We now apply Section \ref{section Covering Types} to the case when $\mathcal{G}$ is the connected classical group $G^0_W$ and 
$\mathcal M = M^0 = \tilde{G}\times G^0$ is a maximal Levi subgroup of $G^0_W$, where $\tilde{G} = \tilde{G}_{H_-}$ and $G = G_V$ as in Section \ref{section Classical groups}. Suppose that $M^0$ is contained in a parabolic subgroup $P$ of $G_W^0$ with unipotent radical $N$. We recall from \cite{Stevens-Miya} on the construction of a covering type in $G^0_W$ of a given cuspidal type in $M^0$. 

Suppose that $\tilde\pi$ is a depth zero supercuspidal representation of $\tilde G$ of the form $\cInd_{ \tilde{ \mathcal J}}^{\tilde G}\tilde{\rho}_{\tilde\xi}$ as in (\ref{depth zero supercuspidal compact induction}). Recall that the inducing type $\tilde{\rho}_{\tilde\xi}$ is constructed  from a cuspidal representation $\tilde{\rho}$ of $\tilde{\mathsf{G}}:=\tilde{\mathcal J}_0/ \tilde{\mathcal  J}_{0+}\cong \GL(m,\mathbb F_q)$ and a character $\tilde\xi$ of $F^\times$. If $\tilde\rho$ is self-dual, i.e., $\tilde{\rho}\cong \tilde{\rho}^\sigma:={}^c\tilde\rho^{\vee}$, then it admits two self-dual extensions $\tilde{\rho}_{\tilde\xi}$ to $ \mathcal J = F^\times  \mathcal J_0$, determined by whether $\tilde{\rho}_{\tilde\xi}(\varpi) = \tilde\xi(\varpi)=\pm 1$. Here we fix $\tilde\pi$ that is induced from $\tilde{\rho}_{\tilde\xi}$ with 
$$\tilde\xi(\varpi)= 1.$$

The construction of covering types is on the level of connected groups. Suppose that $\pi$ is a depth zero supercuspidal representation of $G^0$ of the form $\cInd_{\mathcal{J}}^{G^0}\rho$ where $\mathcal{J}$ is a maximal compact subgroup of $G^0$ with underlying parahoric subgroup $\mathcal{J}_0$ and pro-unipotent radical $\mathcal{J}_{0+}$, and $\rho$ is inflated from a cuspidal representation of $\mathcal{J}/\mathcal{J}_{0+}$.


Denote $\mathcal{J}_{M^0}:=i_M(\tilde{\mathcal J}_0\times  \mathcal{J})$, a maximal compact subgroup of $M^0$, with underlying parahoric subgroup $\mathcal{J}_{M^0,0}$ and pro-unipotent radical $\mathcal{J}_{M^0,0+}$. As in \cite[Sec 3.8]{BHS}, there are two maximal compact subgroups $\mathcal{J}_{W,w}$, for $w\in \{y,z\}$, in $G_W$ whose respective intersections with $M^0$ are both $\mathcal{J}_{M^0}$ and moreover. If we denote by $\mathcal{J}_{W,w,0}$ its  underlying parahoric subgroup with pro-unipotent radical $\mathcal{J}_{W,w,0+}$, then 
$$ \mathcal{J}_P: = \mathcal{J}_{W,y} \cap \mathcal{J}_{W,z} =(\mathcal{J}_{W,y} \cap P)\mathcal{J}_{W,y,0+} = (\mathcal{J}_{W,z} \cap P_-)\mathcal{J}_{W,z,0+}. $$
Note that  $\mathcal{J}_P  =(\mathcal{J}_{W,y} \cap P)(\mathcal{J}_{W,y,0+}\cap N_-) $, and denote 
$$\mathcal{J}_{P,0} = (\mathcal{J}_{W,y,0} \cap P)(\mathcal{J}_{W,y,0+}\cap N_-)\quad\text{and}\quad \mathcal{J}_{P,0+} = (\mathcal{J}_{W,y}\cap N_+)(\mathcal{J}_{{W,y},0+}\cap P_-).$$ One may write down similar decompositions with $y$ replaced by $z$.


Suppose that $\rho_{M^0} = \tilde{\rho} \times \rho$ is a depth zero cuspidal type of $\mathcal J_{M^0} = \tilde{\mathcal J}_0\times \mathcal{J}$. Via the isomorphism of reductive quotients \begin{equation}
\label{JP-and-JM-have same reductive quotient}
 \mathcal{J}_{P}/\mathcal{J}_{P,0+} \cong  \mathcal{J}_{M^0}/\mathcal{J}_{M^0,0+} \cong \tilde{\mathsf G}\times \mathsf S(\mathsf   G_{y}\times \mathsf G_{z}),
\end{equation}
we let $\rho_{P}$ be the inflation of $\rho_{M^0}$ to $\mathcal{J}_{P}$. By \cite[Th 1.1]{Stevens-Miya}, $(\mathcal J_{P},\rho_{P})$ is a covering type of $(\mathcal J_{M^0},\rho_{M^0})$, i.e., it satisfies conditions (i)-(iii) in Section \ref{section Covering Types}.

As a remark, the quotient $\mathcal{J}_{W,y} / \mathcal{J}_{W,y,0+}$ is of the form $\mathsf S(\mathsf G_{W,y}\times \mathsf G_{z})$, where  $\mathsf G_{W,y}$ is the classical group of the same type of $\mathsf G_y$ and containing  $\mathsf M_y:=\tilde{\mathsf G}\times \mathsf G_{y}$ as a maximal Levi subgroup. Similarly,  $\mathcal{J}_{W,z} / \mathcal{J}_{W,z,0+}\cong \mathsf S(\mathsf G_{y}\times \mathsf G_{W,z})$, with $\mathsf G_{W,z}$ containing  $\mathsf M_z:=\tilde{\mathsf G}\times \mathsf G_{z}$ as a maximal Levi subgroup. The quotient in (\ref{JP-and-JM-have same reductive quotient}) is then a maximal Levi subgroup of both $\mathcal{J}_{W,w} / \mathcal{J}_{W,w,0+}$, where $w\in \{y,z\}$.

\subsection{Hecke algebras and modules}
\label{subsection Hecke algebras and modules}


The Hecke algebra $\mathcal{H}(M^0,\rho_{M^0})$ is commutative by \cite[(5.6)]{BK-cover}. In fact, we have in our case $
\mathcal{H}(M^0,\rho_{M^0})\cong \mathbb{C}[\left<\varpi\right>]$ the group algebra of the infinite cyclic group $\left<\varpi\right>$ over $\mathbb{C}$, whose simple modules are parametrized by the group of unramified characters of $F^\times$. Let $\tilde \pi \times \pi$ be a product of supercuspidal representations of $M^0=\tilde{G}\times G^0$. For $s\in \mathbb{C}$, we denote by $D_s$ the simple right $\mathcal{H}(M^0,\rho_{M^0})$-module 
$$\mathbf{M}_{M^0}(\tilde{\pi}|\det|^s\times\pi)
=\Hom_{\mathcal{J}_{M^0}}(\rho_{M^0},\tilde{\pi}|\det|^s\times\pi).$$ 
This module is 1-dimensional over $\mathbb{C}$. Let $Z\in \mathcal{H}(M^0,\rho_{M^0})$ whose support is $\mathcal{J}_{M^0}(\varpi\times I_{V})\mathcal{J}_{M^0} = (\varpi\times I_{V})\mathcal{J}_{M^0}$. Put $Z(\varpi\times I_V)$ simply as $Z(\varpi)$. The action of $Z$ on $D_s$ is given by \cite[(1.15)]{Bl-Bl-SP4}:
\begin{equation}
\label{action of Z on D}
Z\cdot \phi = |\varpi|^{-s}\tilde{\rho}(\varpi)^{-1} \phi \circ  Z(\varpi) = q^{s} \tilde{\rho}(\varpi)^{-1}\phi \circ  Z(\varpi).
\end{equation}
This is just $a_s\phi$ for some $a_s\in \mathbb{C}$. Note that $a_s$ depends on the normalization of $Z$.

We recall the structure of the Hecke algebra $\mathcal{H}(G^0_W,\rho_{P})$ from \cite{Stevens-Miya}. Let $N_{G^0_W}(\mathfrak{s}_M)$ be the normalizer of the inertial class $\mathfrak{s}_M = [\tilde \pi\times \pi]_M$. There are two cases. If $N_{G^0_W}(\mathfrak{s}_M)=M^0$, then $\mathcal{H}(G^0,\rho_{P}) \cong \mathcal{H}(M^0,\rho_{M^0}) $. This is an uninteresting case which we will discard. If $N_{G^0_W}(\mathfrak{s}_M)\supsetneq M^0$, then indeed $[N_{G^0_W}(\mathfrak{s}_M): M^0]=2$, and $\mathcal{H}(G^0_W,\rho_{P,0})$ is isomorphic to the generic Hecke algebra on an infinite dihedral group \cite[Th 1.2]{Stevens-Miya}. More precisely, it is generated by the images of two injective morphisms
$$\mathcal{H}({\mathcal{J}}_{W,w}, \rho_{M^0}) \hookrightarrow \mathcal{H}(G^0_W,\rho_{P}),\,\qquad w\in \{y,z\}.$$
For a cuspidal type of $M^0$ of arbitrary depth, there is an associated tamely ramified character  of ${\mathcal{J}}_{W,w}$ involved in the above embedding; see \cite[Sec 3.12]{BHS} or \cite[Th 1.1]{tam-very-cusp} for example. For depth zero cuspidal types, the result of \cite{Lust-Stevens} implies that this character is trivial.

We now describe the structure $\mathcal{H}({\mathcal{J}}_{W,w}, \rho_{M^0}) $. We first apply a reduction by an application of \cite{Morris-level-zero},
$$\mathcal{H}({\mathcal{J}}_{W,w}, \rho_{M^0}) \cong \mathcal{H}({\mathcal{J}}_{W,w,0}, \rho_{M^0}^0) ,$$
 where $\rho_{M^0}^0$ is a component of the restriction of $\rho_{M^0}$. For complete detail of this reduction for classical groups, see the discussions in \cite[Sec 6.3, 6.4]{Stevens-Miya}. With our classification of cuspidal types in Section \ref{subsection Compact inductions}, when both $\mathsf G_w$ are even orthogonal, $\rho_{M^0} $ is of the form $\tilde\rho\times \mathsf S\rho$, where $\mathsf S\rho$ is either $\rho(\pm)$ in (\ref{(++)on-S(OtimesO)}) or $\rho_\epsilon$ in (\ref{all-zero-on-S(OtimesO)}). Hence $\rho_{M^0}^0$ is of the form $\tilde\rho\times \rho^0$ where $\rho^0 =\rho^0_y\times \rho^0_z$ is in either (\ref{type-G0yTtimesG0z}) or (\ref{all-zero-on-SO-times-SO)}). The description for other types of $\mathsf G_w$ are just analogous.

 The next reduction is just mod p:
 $$\mathcal{H}({\mathcal{J}}_{W,w,0}, \rho_{M^0}^0) \cong \mathcal{H}(\mathsf{G}^0_{W,w}, \tilde{\rho} \times \rho_{w}^0),$$
i.e.,  we have reduced to describe the structure of Hecke algebras of groups over $\mathbb F_q$. We only concern that $\tilde{\rho}$ is self-dual, then $\mathcal{H}(\mathsf{G}^0_{W,w}, \tilde{\rho} \times \rho_{w}^0)$ is generated by an element $T_w$ which satisfies a quadratic equation, say 
\begin{equation}
\label{quadratic relation by Tw}
T_w*T_w = b_w T_w + c_w \mathbbm 1_w
\end{equation}
for some $b_w, c_w\in \mathbb{C}$ (and here $*$ is the convolution product of the Hecke algebra, with a fixed Haar measure and a unit element $\mathbbm 1_w$, the characteristic map on $\mathsf M^0_w$ defined by $\tilde{\rho} \times \rho_{w,0}$). Indeed we will see later in (\ref{explicit-cy}) and (\ref{explicit-cz}) that $c_w$ is positive, and $b_w $ is known to be the form  $c_w^{1/2}(q^{r_w/2}-q^{-r_w/2})$,
where  $c_w^{1/2}$ is a square root of $c_w$ and for some $r_w\in \mathbb{Z}_{\geq 0}$. For now the sign of $c_w^{1/2}$ is unspecified. Finally, we will normalize both $T_y$ and $T_z$ such that $$T_y*T_z = t_P(Z),$$
where $t_P$ is the injection $\mathcal{H}({{M}}^0,\rho_{{M}^0}) \rightarrow \mathcal{H}({{G}^0_W},\rho_{{P}})$.

The parameters $r_w$ for classical groups are calculated in \cite{Lusztig-chevalley}. Suppose that $\tilde{\rho}$  corresponds to $\tilde s\in \hat{\tilde{\mathsf G}}$ with polynomial $\tilde{P}$ which appears as a factor of the polynomial $P$ corresponding to $\rho$, with multiplicity $N_{\tilde P,w}$ which takes the value listed in Section \ref{subsection Finite classical groups}. The values of $r_w$ are listed in \cite{Lust-Stevens} as follows.
\begin{equation}
\label{Lusztig parameter}
\begin{split}
&\text{$\mathsf G^0_w = \SO(2n_w+1)$,  then  $r_w =2m_+ +1$ when $\tilde{s}=1$, and $r_w =2m_- +1$ when $\tilde{s}=-1$.}
\\
&\text{$\mathsf G_w = \SP(2n_w)$,  then  $r_w =2m_+ +1$ when $\tilde{s}=1$, and $r_w =2m_-$ when $\tilde{s}=-1$}
\\
&\text{$\mathsf G_w^0 = \SO_\pm(2n_w)$, then $r_w =2m_+ $ when $\tilde{s}=1$, and $r_w =2m_- $ when $\tilde{s}=-1$}
\\
&\text{$\mathsf G_w =  \mathrm U(N_w)$, then  $r_w = (m_{\tilde P}+1/2)\deg \tilde P$.
}
\end{split}
\end{equation}
We now recall that $T_w$, as a map from $G^0_W$ to the endomorphism algebra of the representation space $\mathbf{V}_{\rho_P})$ of $\rho_P$, is supported on a single $\mathcal{J}_{P}$-double coset $\mathcal{J}_{P}s_w\mathcal{J}_{P}$ with $s_w\in \mathcal{J}_{W,w}\smallsetminus \mathcal{J}_{W,{w'}}$ where $\{w,w'\}=\{y,z\}$. The operator $T_w(s_w)$ leaves the representation spaces $\mathbf V_{\tilde\rho}\cong \mathbf V_{{}^\sigma \tilde\rho}$ of $\tilde\rho$ and $\mathbf V_{\rho}$ of $\rho$ invariant, so we write $T_w(s_w)= \tilde A_w \times A_w$, where 
$$\tilde A_w = T_w(s_w) |_{\bf V_{\tilde\rho}}\quad\text{and}\quad
A_w = T_w(s_w) |_{\bf V_{\rho}}.$$ 
Here $\tilde A_w $ is a non-zero intertwining operator between $\tilde{\rho}$ and ${}^\sigma\tilde{\rho}$.

We follow \cite[Sec 6.2]{stevens-supercusp} and (with respect to a suitable choice of the Hermitian form $h_V$ on $V$) choose 
\begin{equation}
\label{representative-sy-and-sz}
s_y  = \antidiag(I_{},\dot{\mathbbm p}_y,I_{})\qquad \text{ and }\qquad s_z =\antidiag(\varpi_{}^{-1}I_{}, \dot{\mathbbm{p}}_z,  \varpi_{}^{}I_{}).
\end{equation}
Here $\dot{\mathbbm p}_y = \diag(\mathbbm p_y, I_{V_z})$ and $\dot{\mathbbm p}_z = \diag(I_{V_y}, \mathbbm p_z)$, where $\mathbbm p_w$, for $w\in \{y,z\}$, is $I_{V_w}$ unless $\mathsf G_w$ is orthogonal, in which case we lift $\mathbbm p_w$ from $\mathsf G_w \smallsetminus \mathsf G_w^0 = \mathrm O(\mathsf V_w)\smallsetminus \mathrm{SO}(\mathsf V_w)$ with $\mathbbm p_w^2=1$. 
With these representatives $s_w$, we have 
$$T_y(s_y)^2 = T_z(s_z)^2=1,
$$
and that $\varpi_P := s_y s_z$ is a strongly positive element with respect to $(P,\mathcal{J}_{P})$. Put $\varpi_P$ into $t_P(Z) = T_y * T_z$, we have
\begin{equation}
\label{relation between Ay, Az and Z(not t_P(Z))}
T_y(s_y)\circ T_z(s_z) = t_P(Z)(\varpi_P) = Z(\varpi) \Delta_P(\varpi_P)^{-1/2},  
\end{equation}
where $\Delta_P$ is the modular character of $P$, and indeed  $\Delta_P(\varpi_P)^{} = c_y c_z$. 

\subsection{Reducibility of modules and representations}
\label{subsection Calculations with Hecke algebras}

We denote by $X_s$ the right 
$\mathcal{H}(G^0_W,\rho_{P})$-module
$$\mathbf{M}(\tilde{\pi}|\det|^s\times\pi_{})=\Hom_{\mathcal{J}_{P}}(\rho_{P},\tilde{\pi}|\det|^s\times\pi)\cong \Hom_{\text{Mod-}\mathcal{H}(M^0,\rho_{M^0})}(\mathcal{H}(G^0_W,\rho_{P}),D_s),$$
where the last isomorphism is given by the commutative diagram (\ref{commutative diagram}).

We follow the idea in \cite[Sec 1]{Bl-Bl-SP4} that leads us to a calculation towards our major formula (\ref{traced_difference_Dy_Dz}) in the next section on the difference between certain orthogonality relations of $\tilde\rho\times \rho_y$ and of $\tilde\rho\times \rho_z$. The key step in (\ref{Final comparison of eigenvalues}) below is to compare the eigenvalues $a_s$ of $t_P(Z)$ and those of $T_w$ on the space $X_s$ via (\ref{relation between Ay, Az and Z(not t_P(Z))}). We first recall an important result from \cite[(1.13)]{Bl-Bl-SP4}.

\begin{prop}
\label{Blondel-Blasco about product of eigenvalues}
The module $X_s$ is reducible {if and only if} the product of eigenvalues of $T_y$ and of $T_z$ on $X_s$ is equal to the scalar $a_s$ of $Z$ on $D_s$. 
\qed\end{prop}

\begin{rmk}
\label{remark by Blondel}
The main result of \cite{Bl-Bl-SP4} applies to SP(4), but the above proposition applies also to other classical groups. All we need is that the Hecke algebra $\mathcal{H}(G^0_W,\rho_{P})$ is a free $\mathcal{H}(M^0,\rho_{M^0})$-module of rank 2, so that the right $\mathcal{H}(G^0_W,\rho_{P})$-module induced from a character of $\mathcal{H}(M^0,\rho_{M^0})$ is of rank 2 and contains an invariant line in the case of reducibility, i.e., $\mathcal{H}(G^0_W,\rho_{P})$ acts on the line by scalars. \footnote{Here we thank Corinne Blondel for the explanation.} Moreover, the proposition in \cite{Bl-Bl-SP4} says more: $X_s$ is completely reducible if and only if $a_s$ is a product of eigenvalues in two different ways. We do not need this statement for our purpose.
\qed\end{rmk}

By this proposition, and with (\ref{action of Z on D}), (\ref{relation between Ay, Az and Z(not t_P(Z))}), we obtain a relation
$$
q^s \tilde{\rho}(\varpi)^{-1}\Delta_P(\varpi_P)^{1/2}  T_y(s_y)\circ T_z(s_z) = \text{ product of eigenvalues }\lambda_y\lambda_z\text{ of $T_y$ and $T_z$}.
$$
By restricting to $\mathbf V_{\tilde\rho}$, the relation implies that
\begin{equation}
\label{Final comparison of eigenvalues}
q_E^s \tilde A_y = \lambda_y\lambda_z \Delta_P(\varpi_P)^{-1/2} \tilde{\rho}(\varpi) \tilde A_z.
\end{equation}
The two eigenvalues of $T_w$ are of the form $c_w^{1/2}q^{r_w}$ or $-c_w^{1/2}q^{-r_w}$, so that the product of eigenvalues is of the form 
$$c_y^{1/2}c_z^{1/2}q^{\pm (r_y+r_z)}
\qquad\text{or}\qquad 
-c_y^{1/2}c_z^{1/2}q^{\pm (r_y-r_z)}.$$
The relation (\ref{Final comparison of eigenvalues}) is the key step relating points of reducibility and parameters of Hecke algebras, resulting in Corollary \ref{cor the points of reducibility} below.

To state this corollary, let's first define $\nu\in \{\pm 1\}$ as follows. If $b_w\neq 0$, we  call $b_w  \tilde A_w$ the {\bf canonical normalization} of the intertwining operator $\tilde A_w$, and denote 
\begin{equation}
\label{introducing dw}
d_w :=b_w  \tr \tilde A_w \qquad
\text{for }w\in \{y,z\}.
\end{equation}
If both $d_w\neq 0$, we denote
$$\nu = \sgn(d_y )\sgn(d_z )\tilde{\rho}(\varpi). $$
When either $r_y$ or $r_z$ is 0, we can take $\nu$ to be either $\pm 1$. From (\ref{Final comparison of eigenvalues}), we conclude the following reducibility result about the parabolic induction
$I(s,\tilde{\pi},\pi)=\iota_{P}^{G^0_W}(\tilde{\pi}_{}|\det|^s\times \pi)$, with $s\in \mathbb{C}$.
\begin{cor}
\label{cor the points of reducibility}
The points of reducibility of $I(s,\tilde{\pi},\pi)$ in the domain (\ref{fundamental domain of reducibility}) are
$$ \pm \frac{r_y+ \nu  r_z}{2}\qquad\text{and}\qquad
\pm \frac{r_y- \nu  r_z}{2} + \frac{\pi \sqrt{-1}}{m\log q}.$$
(Note here $\tilde\pi$ is a representation of $\GL(m)$.)
\end{cor}

Denote by $\tilde{\pi}'$ the self-dual representation differs from $\tilde{\pi}$ by twisting by an unramified character. The  points of reducibility of $I(s,\tilde{\pi}',\pi)$ are the same as in Corollary \ref{cor the points of reducibility}, but with $ \nu  $ switched into $- \nu  $.

\proof (of Corollary \ref{cor the points of reducibility}) 
The real parts of the points are already given by \cite[Prop 3.12]{Blondel-Weil}, which are 
$$\pm \frac{r_y+  r_z}{2}\qquad\text{and}\qquad
\pm \frac{r_y- r_z}{2}$$
as a multi-set of four numbers. When either $r_y$ or $r_z=0$, then these real parts are all the same up to a sign, and the corollary is obvious. When both $r_y,r_z\neq 0$, the exact points of reducibility can be computed from (\ref{Final comparison of eigenvalues}). \qed

Note that $b_w=0$ if and only if $r_w=0$. Using the list (\ref{Lusztig parameter})), this happens only when 
\begin{equation}
\label{necessary conditions for b neq 0}
\begin{split}
\text{either }&(\mathsf G_w\text{ is even orthogonal, }\tilde\rho = 1\text{, and }m_{+,w}=0)
\\
\text{or }&(\mathsf G_w\text{ is symplectic or even orthogonal, }\tilde\rho = \omega\text{, and } m_{-,w}=0).
\end{split}
\end{equation}
In contrast, $r_w\neq 0$ implies that $m_{\pm ,w}\neq 0$ for $\rho$ as in (\ref{necessary conditions for b neq 0}), which further implies that $\rho_{M^0}$ extends to $\rho_{M}$.

\subsection{More explicit calculations}
\label{subsection Proof of the main result}

We show that, for $w\in \{y,z\}$, the canonical normalization $d_w$ is a sum over, and depending on, only character values of $\tilde\rho\times \rho_w$. We will also determine the sign $\nu$ as a consequence.

We use a direct computation in \cite[(1.3)]{Bl-Bl-SP4} for the coefficients $b_w$ and $c_w$ of the quadratic relation (\ref{quadratic relation by Tw}) for $T_w$. Let $s_w$ be defined as in (\ref{representative-sy-and-sz}). With the Haar measure chosen such that the volume of $\mathcal{J}_{P}$ is 1, we have 
 \begin{equation}
 \label{explicit-cy}
c_y =[\mathcal{J}_P^+:s_y\mathcal{J}_P^- s_y],\qquad \text{where }\mathcal{J}^+_P = \mathcal{J}_{P}\cap N\text{ and }\mathcal{J}^-_P=\mathcal{J}_{P}\cap N_-,
\end{equation}
and if we decompose $j\in s_y\mathcal{J}_P^+s_y$ into an element of the form $\mathcal{J}_{P}s_y\mathcal{J}_{P}$, as 
$$j = u_1 s_y m(j) u_2, \qquad u_1,u_2\in \mathcal{J}_P^+,\,m_y(j)\in M,$$
then 
\begin{equation}
\label{by-end}
b_y(\tilde A_y\times I_V)=\sum_{j\in \frac{s_y\mathcal{J}_P^+s_y\cap \mathcal{J}_Ps_y\mathcal{J}_P}{\mathcal{J}_P^- }} \rho_{M} (m_y(j)).
\end{equation}
where $\rho_{M} = \tilde\rho \times \rho_y\times \rho_z $ as a cuspidal type extending $\rho_{M^0}$. With the similar calculation as in \cite[Sec 4.1]{Blondel-Tam-2020}, the ranging set in the sum (\ref{by-end}) is 
\begin{equation}
\label{Sy-Set}
\mathcal{S}_y :=\left\{ (X_y,Y)\in \Hom_{\mathbb F_q}(\mathsf V_y, \mathsf H_+) \times \Hom_{\mathbb F_q}(\mathsf H_-, \mathsf H_+)^\times  \text{ satisfying  
(\ref{XalphaX=Y-alphaY}) mod } \mathcal{P}_F\right\}, 
\end{equation}
where $\Hom_{\mathbb F_q}(\mathsf H_-, \mathsf H_+)^\times$ means the subset of invertible maps in $\Hom_{\mathbb F_q}(\mathsf H_-, \mathsf H_+)$, and is simply identified with $\tilde{\mathsf G} = \GL(m)$ with $m=\dim\mathsf  H_-$. Taking trace on (\ref{by-end}), we have 
\begin{equation}
\label{Sy-Set-orthogonality-relation}
b_y\tr \tilde A_y \deg  \rho_y =\sum_{(X_y,Y)\in \mathcal{S}_y}
\tr  \tilde{\rho}(Y) . \tr {\rho}_y(I_{\mathsf V_y}-{}^\alpha X_yY^{-1}X_y).
\end{equation}
Similarly, we have
\begin{equation}
 \label{explicit-cz}
c_z = [s_z \mathcal{J}_P^-s_z^{-1}:\mathcal{J}_P^+],
\end{equation}
and if we write $j = u_1  m_z(j) s_z u_2$ for some $ u_1,u_2\in \mathcal{J}_P^-$ and $m_z(j)\in M$, then
\begin{equation}
b_z(\tilde A_z \times I_V)=\sum_{j\in \frac{s_z^{-1}\mathcal{J}_P^-s_z\cap \mathcal{J}_Ps_z\mathcal{J}_P}{\mathcal{J}_P^+ }} \rho_M (m_z(j)) 
\end{equation}
which implies that 
\begin{equation}
\label{Sz-Set-orthogonality-relation}
b_z\tr \tilde A_z\deg \rho_z =\sum_{(X_z,Y)\in \mathcal{S}_z}\tr \tilde{\rho}(Y).\tr {\rho}_z(I_{\mathsf V_z}-{}^\alpha X_zY^{-1}X_z)
\end{equation}
 where 
$\mathcal{S}_z$ is defined similarly as in 
(\ref{Sy-Set}), with $\mathsf H_+$ and $\mathsf H_-$ exchanged and $y$ replaced by $z$.

Let's look at the difference between (\ref{Sy-Set-orthogonality-relation}) and (\ref{Sz-Set-orthogonality-relation}),
\begin{equation}
\label{traced_difference_Dy_Dz}
\begin{split}
&b_y\tr \tilde A_y \deg \rho_y- b_z\tr \tilde A_z\deg \rho_z   = 
 \sum_{Y\in \mu_F}
\tr \tilde{\rho}(Y)\cdot
\\
&\Bigg(\sum_{
\begin{smallmatrix}
X_y\text{ such that}
\\
(X_y,Y)\in\mathcal S_{y}
\end{smallmatrix}
}
\tr {\rho}_y(I_{\mathsf V_y}-{}^\alpha X_yY^{-1}X_y) -
\sum_{
\begin{smallmatrix}
X_z\text{ such that}
\\
(X_z,Y)\in \mathcal S_{z}
\end{smallmatrix}
} 
 \tr {\rho}_z(I_{\mathsf V_z}-{}^\alpha X_zY^{-1}X_z)\Bigg).
\end{split}
\end{equation}
Note that $\tilde{A_y}$ and $\tilde{\rho}(\varpi)\tilde A_z$ differ by a sign; if both $d_w\neq 0$, this sign is $\nu \sgn(b_y b_z)$. Also, we have fixed $\tilde\rho(\varpi)=1$. In the typically almost symmetric case, i.e., if $\rho_y$ is a companion of $\rho_z$, we can simplify (\ref{traced_difference_Dy_Dz}) a bit as \begin{equation}
\label{cancellation in traced_difference_Dy_Dz}
\begin{split}
&\sgn(d_y) |\tr\tilde A_y| \deg \rho _y( |b_y| - \nu | b_z|)   = 
 \sum_{Y\in \mu_F}
\tr \tilde{\rho}(Y)\cdot
\\
&\Bigg(\sum_{
\begin{smallmatrix}
X_y\text{ such that}
\\
(X_y,Y)\in\mathcal S_{y}\text{ and }
X_y\neq 0
\end{smallmatrix}
}
\tr {\rho}_y(I_{\mathsf V_y}-{}^\alpha X_yY^{-1}X_y) -
\sum_{
\begin{smallmatrix}
X_z\text{ such that}
\\
(X_z,Y)\in \mathcal S_{z}\text{ and }X_z\neq 0
\end{smallmatrix}
} 
 \tr {\rho}_z(I_{\mathsf V_z}-{}^\alpha X_zY^{-1}X_z)\Bigg)
\end{split}
\end{equation}
using Lemma \ref{properties as in Waldspurger book} .

We need a few properties about elements of the form $I - {}^\alpha X Y^{-1}X$ in the formula above when $\deg\tilde \rho=1$. For an orthogonal or a symplectic group $\mathsf G = \mathsf G(\mathsf V)$, identify $\Hom_{\mathbb F_q}(\mathsf H_-, \mathsf H_+)^\times$ with ${\mathbb F_q}^\times$ and denote $\Hom_{\mathbb F_q}(\mathsf V, \mathsf H_+) $ by $\mathsf V^*$ so that $\alpha(\mathsf V^*) = \mathsf V$. 

\begin{lem}
\label{lemma of similitudes}
Consider elements $I - {}^\alpha X Y^{-1}X\in \mathsf G(\mathsf V)$, where 
$$(X,Y)\in \mathcal{S} :=\left\{ (X,Y)\in( \mathsf V^* \smallsetminus\{0\}) \times \mathbb F_q^{\times }\text{ satisfying  
(\ref{XalphaX=Y-alphaY})} \right\}.$$ 
\begin{enumerate}[(i)]

\item We have $\det (I - {}^\alpha X Y^{-1}X) = -\epsilon_{\mathsf G}$.
\label{determinant of conjugacy class is 1 or -1}

\item The elements $I - {}^\alpha X Y^{-1}X$, with $(X,Y)$ ranging over $\mathcal S$, lie in only two distinct $\mathsf G$-conjugacy classes $\mathcal C_{{\square }}$ and $\mathcal  C_{{\centernot\square }}$, where
$$\mathcal  C_{{\square \text{ (resp. }\centernot\square)} }:=\{g (X_{\square \text{ (resp. }\centernot\square)})g^{-1}:g\in \mathsf G\}$$ for two elements $X_{\square}$ and $X_{\centernot\square}$ in $\mathsf G$.
\label{two conjugacy classes}

\item That  $I - {}^\alpha X Y^{-1}X$ lies in which $\mathsf G$-conjugacy class depends on whether $Y\in \mathbb F_q^{\times 2}$ or not.
\label{square or not}

\item Let $h\in \mathsf{GG}(\mathsf V)$ be the element in the similitude group of $\mathsf{G}$ defined in Lemma \ref{properties as in Waldspurger book}, then   $\Ad(h)$  exchanges  $\mathcal C_{{\square  }}$ and $\mathcal  C_{{\centernot\square} }$.
\label{definition of h}

\end{enumerate}
\end{lem}
\proof
For (\ref{determinant of conjugacy class is 1 or -1}), we have $\det (I - {}^\alpha X Y^{-1}X) = 1 - {}X^\alpha X Y^{-1} $, and we use (\ref{XalphaX=Y-alphaY}) to show that it is equal to $ 1 - (Y - {}^\alpha Y) Y^{-1} = -\epsilon_{\mathsf G}$. For (\ref{two conjugacy classes})-(\ref{definition of h}), we prove all these together using the following setup. Consider the action of $\mathsf{GM}:=\tilde{ \mathsf G}(\mathsf H_-)\times \mathsf{GG(V)}$ on the set 
$ \mathcal S$  by 
\begin{equation*}
(x,g )\cdot (X,Y) = (x Xg^{-1}, x Y {}^\sigma x^{-1} \lambda(g)^{-1} ).
\end{equation*}
The map
$$\mathcal{N}:\mathcal S\rightarrow \mathsf{ G},\,\quad (X,Y)\mapsto I-{}^\alpha X Y^{-1}X$$ 
is constant along $(\tilde{ \mathsf G}(\mathsf H_-)\times 1)$-orbits, and the composition mapping $(X,Y)$ into the $\mathsf{G}$-conjugacy class of $\mathcal{N}(X,Y)$ is constant along $\mathsf{GM}$-orbits. We now assume that $\dim \mathsf H_-=1$, then $ {}^\sigma x = x^{-1} $ for all $x\in \tilde{ \mathsf G}(\mathsf H_-) = \mathsf{GL}(1)$. For $g\in \mathsf G$, we have $\Ad(g)(\mathcal N ( X,Y) ) =\mathcal N ( (1,g)\cdot(X,Y) ) $. Since the image of $\lambda|_{\mathsf G}$ is $\mathbb{F}_q^{\times2}$, we see that $\mathsf G$-conjugacy leaves $Y\bmod \mathbb{F}_q^{\times2}$ invariant. If we take $g=h \in \mathsf{ GG(V)}$, then we see that $\Ad(h)$ brings  $\mathcal N(X,Y) $ to another $\mathsf G$-conjugacy class.

To show that $\mathcal N(\mathcal S)$ only covers two $\mathsf G$-conjugacy classes, first consider when $\mathsf G$ is orthogonal, and note that $X{}^\alpha X\neq 0$ (since  $1 - {}X^\alpha X Y^{-1}=-1$). Suppose that  $(X',Y')\in \mathcal{S}$ such that $(X{}^\alpha X) (X'{}^\alpha X')^{-1}=x^2$ for some $x\in \mathbb F_q^\times$. By Witt's Theorem, the map $X\mapsto xX'$ can be extended to an isometry $\mathsf T:\mathsf V^*\rightarrow \mathsf V^* $, and there is $g\in \mathsf {G(V)}$ such that $\mathsf T(v) = v g^{-1}$ for all $v\in  \mathsf V^* $. In particular, we have $ X' =  x Xg^{-1}$, and hence $(x,g)\cdot (X,Y) = (x Xg^{-1}, x^2 Y )$.  By taking the map $\mathcal N$ to this equality, we obtain (\ref{two conjugacy classes}) and (\ref{square or not}). If $(X{}^\alpha X) (X'{}^\alpha X')^{-1}\notin \mathbb{F}_q^{\times2}$, then we replace $(X',Y')$ by $(1,h)\cdot (X',Y')$. Arguing as before, we obtain (\ref{definition of h}) this time.

If $\mathsf G$ is symplectic, we first take $Z\in \mathsf V^*$ such that $Z{}^\alpha X = 1$. Suppose that  $(X',Y')\in \mathcal{S}$, and take $Z'\in \mathsf V^*$ such that $Z'{}^\alpha X' = 1$. If $Y /Y' = x^2$ for some $x\in \mathbb F_q^\times$, then the map $X\mapsto xX',Z\mapsto x^{-1}Z'$ is an isometry from the subspace of $\mathsf V^*$ spanned by $\{X,Z\}$ onto the subspace spanned by $\{X',Z'\}$, and it extends to an isometry $g\in \mathsf G$ such that $(x,g)\cdot (X',Y') = (X,Y)$. All the statements can then be proved by arguing as in the orthogonal case. \qed

\subsection{Canonical normalization of intertwining operators}
\label{subsection Canonical normalization of intertwining operators}

We have seen, in (\ref{Sy-Set-orthogonality-relation}) or (\ref{Sz-Set-orthogonality-relation}) above, a calculation unfolding the character of the canonical normalization $d_w$ (\ref{introducing dw}) as a sum involving only character values of $\tilde\rho$ and $\rho$. It happens that this quantity is actually $\neq 0$ as long as $b_w\neq 0$, by the Lemma below.

\begin{lem}
$d_w\neq 0$ if and only if $b_w\neq 0$ (which in turn is equivalent to $r_w \neq 0$).
\end{lem}
\proof If $\deg\tilde \rho = 1$, then $\tilde A_w$ is just a non-zero scalar, and the Lemma is obvious. We henceforth assume that $\tilde\rho$ is fixed with $m:=\deg\tilde \rho > 1$ and $b_w\neq 0$. We first apply a simple reduction: put $\tilde A_w = \tilde A_w(\rho) $ to emphasize the dependence of $\tilde A_w$ on $\rho$, and denote $\tilde A_w(0)$ when $\dim V=0$. As a non-zero intertwining operator between $\tilde\rho$ and ${}^\sigma\tilde \rho$, we have that  $ \tilde A_w(\rho_0)$ is a non-zero scalar multiple of $\tilde A_w(0)$. Hence it is enough to show that $d_w \neq 0$ with $\tilde A_w = \tilde A_w(0)$. 

We emphasize the setup: viewing $\tilde\rho$ as a cuspidal type of a depth zero supercuspidal representation of $\tilde{G} = \GL(m)$ and $\tilde{G}$ as a maximal Levi subgroup of $G_W$, if $m$ is even, we assume $G_W$ is orthogonal or symplectic, while if $m$ is odd we assume $G_W$ is unitray.

In the symplectic or unitary case, to show that $d_w\neq 0$, we simply reiterate the proof in \cite[Prop 9.9(i)]{MR-classical}. We recall briefly the idea here as we need it for the orthogonal case. To begin with, we pick a skew-symmetric (resp. skew-Hermitian) $K\in \tilde{\mathsf G} = \GL(m)$ such that $\{YK^{-1}:{}^\alpha Y = Y\}$ is precisely the subset of skew-symmetric (resp. skew-Hermitian) matrices in $\tilde{\mathsf G}$, which form a unique equivalence class. By (\ref{Sy-Set-orthogonality-relation}),
$$d_w = \sum_{{}^\alpha Y = Y}\tilde{\rho}(Y)
=
\sum_{S\in\mathsf {\tilde G}/ G(K)}\tilde{\rho}(SK^{-1} {}^tS K)$$ where $G(K)$ is the classical group over $\mathbb F_q$ defined by $K$. If we denote the action $\sigma_K: S\mapsto K^{-1} {}^tS^{-1} K$, which is an involution by choosing a suitable $K$, then 
\begin{equation}
\label{twisted sum in Murnaghan-Repka}
\sum_{S\in\mathsf {\tilde G}/ G(K)}\tilde{\rho}(SK^{-1} {}^tS K) = |G(K)|^{-1}\sum_{S\in\mathsf {\tilde G}}\tilde{\rho}(SK^{-1} {}^tS K) = |G(K)|^{-1}\sum_{S\in\mathsf {\tilde G}}\tilde{\rho}(S{}^{\sigma_K} S^{-1}).
\end{equation}
By Schur's Lemma, the last operator is a multiple of $\tilde{\rho}(\sigma_K)$, an extension of $\tilde{\rho}$ to $\tilde{\mathsf G}\rtimes \left<\sigma_K\right>$. Since $\tilde{\rho}(\sigma_K)$ is just a non-zero multiple $\tilde{\rho}(\sigma)$ (as intertwining operators between $\tilde \rho$ and ${}^\sigma \tilde\rho$), one can show that the trace of (\ref{twisted sum in Murnaghan-Repka}) 
is a non-zero multiple of $(\tr\tilde{\rho}(\sigma))^2$. We then apply the twisted Deligne-Lusztig character formula in \cite{DM-non-conn} to show that $\tr\tilde{\rho}(\sigma)$ is a non-zero multiple of $\tilde\xi(\sigma)$, where $\tilde\xi$ is the self-dual character of an elliptic torus $\mathsf T $ that induces to $\tilde\rho$. Knowing that $\tilde\xi(\sigma)\neq 0$, this completes the proof in the symplectic and odd unitary cases.

In the orthogonal case, we take $K = \antidiag(1,\dots,1)$, and $YK^{-1}$ is symmetric. The difference here from the previous cases is that there are two congruence classes of symmetric matrices in $\tilde{\mathsf G}$, which can be represented by $K$ and $ZK$, where $Z = \diag(1_{n-1},\zeta,1_n)$ with $\zeta\in \mathbb F_q^\times \smallsetminus \mathbb F_q^{\times2}$. We hence have
$$\sum_{{}^\alpha Y = Y}\tilde{\rho}(Y) = \sum_{S\in\mathsf {\tilde G}/ G(K)}\tilde{\rho}(SK^{-1} {}^tS K) +\sum_{S\in \mathsf {\tilde G}/G(ZK)} \tilde{\rho}(SZK^{-1} {}^tS K),$$
and so 
$$d_w = |G(K)|^{-1}\sum_{S\in\mathsf {\tilde G}}\tr\tilde{\rho}(SK^{-1} {}^tS K) +|G(ZK)|^{-1}\sum_{S\in \mathsf {\tilde G}}\tr \tilde{\rho}(SZK^{-1} {}^tS K)$$
The first sum is non-zero by using arguments similar to the symplectic and unitary cases. The second sum is 0, because $Z\sigma $ can never be $\mathsf {\tilde G}$-conjugate to any element in $\mathsf T \sigma$.
\qed

\section{The final arguments}
\label{subsection The final argument}

While in general it is difficult to compute individual character values involved on the right side of (\ref{cancellation in traced_difference_Dy_Dz}), our typically almost symmetric condition defined in (\ref{condition_rhoy_equals_rhoz}), together with the expected number of supercuspidal representations in a given packet provided in (\ref{number of supercuspidals in a packet}) and (\ref{number of supercuspidals in a packet even SO}), greatly simplify the situation and allow us to easily prove our main results in Section \ref{subsection The main result}.

In the situation of Proposition \ref{main result for unramified unitary}, take $\tilde{\rho}$ such that $\tilde{\rho}|_{\tilde{\mathcal{J}}_0}$ corresponds to a self-dual polynomial $\tilde P$ with degree $>1$, and take $\rho_y = \rho_z$ whose corresponding polynomial $P$ has a factor $\tilde P$. On the RHS of (\ref{cancellation in traced_difference_Dy_Dz}), we see that the difference in the big brackets is 0; while on the LHS, since $|b_y| = |b_z|\neq 0$, we conclude that $\nu=1$. Therefore, $I(s,\tilde{\pi},\pi)$ reduces at $s= m_{\tilde P} + 1/2$ when $\tilde{\rho}(\varpi)=1$. i.e., we have proved Proposition \ref{main results on reducibility points}(\ref{main results on reducibility points even unitary}).

On the parameter side, the above result determines the parity (whether it is orthogonal or symplectic) of a parameter for $\GL(n)$. Indeed, with $\tilde\varphi = \Ind_{\mathcal W_E}^{\mathcal W_F}\tilde{\xi}$ and $\tilde{\xi}(\varpi)=1$, we have $a_{\tilde\varphi} = 2m_{\tilde P}$. The parity of $\tilde\varphi$ is hence determined: it is symplectic when $m$ is even ($G$ is symplectic or orthogonal), and is orthogonal when $m$ is odd ($G$ is unitary). This proves Proposition \ref{main result for unramified unitary}.

When $G$ is symplectic, we assume that $m_-\neq 0$ (or otherwise $c=d=-1$) and put $\rho_y = \rho_z = \rho(m_+,\pm m_-,(m_{\tilde P})_{\tilde P})$. The RHS of (\ref{traced_difference_Dy_Dz}) is 0, and so $\nu=1$, when $|b_y|=|b_z|\neq 0$. This implies that $I(s,\omega_1,\pi)$ reduces at $s= 2m_{-} $, and so $c = 4m_--1$ by (\ref{SL2-component-and-reducibility}), i.e., we have proved Proposition \ref{main results on reducibility points}(\ref{main results on reducibility points symplectic}). By assuming the expected number given by (\ref{number of supercuspidals in a packet}) of supercuspidal representations in the packet, we have found the two (and one when $c=-1$) supercuspidal representations in the packet $\Pi_{[[a,-1,c,-1,(2m_{\tilde P},-1)_{\tilde P}]]}$. The other two supercuspidal representations, given by $\{\rho_y , \rho_z \}  = \{\rho(m_+, m_-,(m_{\tilde P})_{\tilde P}),\allowbreak  \rho(m_+, - m_-,(m_{\tilde P})_{\tilde P})\}$, must be lying in the companion packet $\Pi_{[[a,-1,-1,c,(2m_{\tilde P},-1)_{\tilde P}]]}$. This proves Proposition \ref{main-result-symplectic}.

For even orthogonal groups, the proof of the first statement in Proposition \ref{main-result-orthogonal} is similar to Proposition \ref{main-result-symplectic}, i.e., we have found the four (the expected number in (\ref{number of supercuspidals in a packet}) supercuspidal representations in the packet $\Pi_{[[a,-1,c,-1,(2m_{\tilde P},-1)_{\tilde P}]]}$ given by $\rho_y = \rho_z$. For the second statement in Proposition \ref{main-result-orthogonal}, we take $\rho_y= \det\otimes \rho_z$. We now note that 
$$\rho_w(I - {}^\alpha X_w Y^{-1}X_w)\neq 0,\quad \text{for all } (X,Y)\in \mathcal S_w,\,w\in \{y,z\},$$
 since otherwise the LHS of (\ref{traced_difference_Dy_Dz}) would become zero, and so $\pi_\rho$ with $\rho_y= \det\otimes \rho_z$ would also belong to $\Pi_{[[a,-1,c,-1,(2m_{\tilde P},-1)_{\tilde P}]]}$, which is impossible because we have exhausted all supercuspidals in this packet. These arguments imply that $\nu=-1$ for both $\tilde\rho =1$ and $\tilde\rho = \omega$, and hence the second statement of Proposition \ref{main-result-orthogonal}.

For the third statement, we use the statement in Lemma \ref{properties as in Waldspurger book} that 
$\rho(m_+,m_-,(m_{\tilde P})_{\tilde P})\circ \Ad(h) \cong \rho(m_+,-m_-,(m_{\tilde P})_{\tilde P})$,
where $h\in \mathsf{GO}(\mathsf V_w)$ the group of orthogonal similitudes with factor $\lambda$ such that $\omega\circ \lambda(h)=-1$, and that $\Ad(h)$ exchanges 
the $\mathsf G$-conjugacy classes $\mathcal C_{{\square  }}$ and $\mathcal  C_{{\centernot\square} }$ by Lemma 
\ref{lemma of similitudes}. If we take $\rho_y = \rho_z\circ \Ad(h)$, then (\ref{traced_difference_Dy_Dz}) implies that
$$\text{$\nu=1$ when $\tilde\rho =1$,\quad
 and \quad $\nu=-1$ when $\tilde\rho = \omega$.}
 $$
This proves the third statement. The last statement is then similar, by taking $\rho_y = \det\otimes (\rho_z\circ \Ad(h))$. This completes the proof of Proposition \ref{main-result-orthogonal}.

The proof of Proposition \ref{main-result-orthogonal-odd-odd} is very similar to Proposition \ref{main-result-orthogonal}. We skip the details.

We remark that when $|b_z|=0$ (or if $|b_y|=0$ then we exchange $y$ and $z$), i.e., equivalently $r_z=0$, the real parts of the reducibility points are all equal up to a sign. Using relation (\ref{SL2-component-and-reducibility}), we conclude that $a_{\tilde\varphi} = b_{\tilde\varphi}$, i.e., we don't have to distinguish $\tilde \varphi$ and $\tilde \varphi'$ in the component $\tilde \varphi[[a_{\tilde \varphi}]]\oplus \tilde \varphi'[[b_{\tilde \varphi}]]$ of $\varphi$.

\subsection{Comparison with the LLC for GL(n)}
 \label{(Comparison with the LLC for GL(n))}

In this subsection, we discuss how Proposition \ref{main result for unramified unitary} can be easily deduced by accepting the LLC for $\tilde G = GL(m)$. Recall that the parity of a supercuspidal representation of a general linear group is the same as its parameter, and can be defined intrinsically using respectively symmetric, exterior, or Asai L-functions using the result in \cite{Hen-ext-sym}. 

Let's compare our parity result with the explicit LLC for $\tilde G $ in \cite[Sec 3]{BH-ET3}. Suppose that $\tilde\pi$ is a depth zero supercuspidal representation of $\tilde G$ and is self-dual: $\tilde\pi \cong \tilde\pi^\vee$. If $m>1$, then $m$ is even. Furthermore, if $\tilde\pi$ is compactly induced from a cuspidal type $\tilde{\rho}_{\tilde\xi}$ as in (\ref{depth zero supercuspidal compact induction}), then by \emph{loc. cit.} its parameter is $
\tilde\varphi \cong \Ind_{\mathcal W_E}^{\mathcal W_F}\tilde\xi\omega_0^E$, where $E/F$ is an unramified extension of degree $m$ and $\omega_0^E$ is the quadratic unramified character of $E^\times$. If $\tilde\xi(\varpi)=1$, i.e. $\tilde\xi$ is conjugate-orthogonal, which implies that $\tilde\xi\omega_0^E$ and also $
\tilde\varphi$ are conjugate-symplectic. By \cite{Shahidi-twisted-endos} (or see the summary in \cite[Sec 4.3]{Shahidi-Cogdell-2014}), this implies that the exterior L-function $L(s,\tilde\varphi, \wedge^2) = L(s,\tilde\pi, \wedge^2)$ has a pole at $s=0$. By \cite{Moeg-exhaustion}, if $\tilde\varphi[[a]]$ is a component of a self-dual parameter, then $a$ is even. 

Parallel arguments hold when we have a quadratic unramified extension $F/\Fo$ and  $\tilde\pi $ is conjugate self-dual: $\tilde\pi \cong {}^c\tilde\pi^\vee$ where $1\neq c\in \Gal(F/\Fo)$. In this case, $m=[E:F]$ is odd and $
\tilde\varphi \cong \Ind_{\mathcal W_E}^{\mathcal W_F}\tilde\xi$, i.e., $\omega_0^E$ is trivial. If $\tilde\xi(\varpi)=1$, then by \cite[Prop 4.31]{Shahidi-Cogdell-2014} the $+$-Asai L-function $L(s,\tilde\varphi, \mathrm{Asai}_+) = L(s,\tilde\pi, \mathrm{Asai}_+)$ has a pole at $s=0$. By \cite{Moeg-base-change}, $a$ is even (resp. odd) if $\tilde\varphi[[a]]$ is a component of a conjugate self-dual parameter of an even (resp. odd) unitary group.

\subsection{Example: Cuspidal quadratic-unipotent representations of Sp(4)}
\label{subsection Examples of quadratic-unipotent representations}

Consider the following parameters for $G=\SP(4)$, 
$$\varphi_i = 1\oplus \omega_i \otimes \St{[[3]]} = [[1,-1,3,-1]]\quad(\text{resp.}\quad [[1,-1,-1,3]]) ,$$
where $i=1 $ (resp. $i=2$). Each packet $\Pi_{\varphi_i}$ contains two supercuspidal representations. In \cite[e.g. 9.4]{Lust-Stevens}, we have distinguished the 4 supercuspidal representations in the union $\Pi_{\varphi_1}\sqcup \Pi_{\varphi_2}$. We will show below which two representations belong to which $\Pi_{\varphi_i}$.

Recall that each representation is compactly induced from a cuspidal representation $\rho_y\times \rho_z$ of $\mathcal J $ such that $\mathcal J/\mathcal J_{0+}\cong \SP(2)\times  \SP(2)$, where each $\rho_w$, for $w\in \{y,z\}$, is chosen from the two cuspidal representations $\rho_\pm = \rho(0,\pm 1)$, both of degree $(q-1)/2$.

 We take $\tilde\rho = \omega$ the quadratic character of $\mathbb F_q^\times$, then 
 $$|b_y| = |b_z| = q^{1/2}(q^2-1).$$
Also, if $\rho_y = \rho_z$, then the RHS of (\ref{traced_difference_Dy_Dz}) is 0, and so $\nu=1$, which means that 
$$ \pi\in \Pi_{\varphi_1}
\quad\text{ if }\quad
\rho_y = \rho_z.$$
Consequently, $\pi\in \Pi_{\varphi_2}
$ if $
\rho_y \neq \rho_z$.

To verify the last statement by directly computing that $\nu=-1$, we have to invoke from \cite[Table 5.4, p.58]{Bonnafe-SL2} the character values of $\rho_{\pm}$ on elements in Sp(2) of the form 
$$I - {}^\alpha X Y^{-1}X,\quad \text{where}\quad
X{}^\alpha X = Y- {}^\alpha Y.$$
A direct computation shows that these elements lie in only two $\SP(2)$-conjugacy classes
 \begin{equation*}
 \begin{split}
 &\mathcal C_\square\qquad \text{ represented by } \begin{bmatrix}
1 & 1 \\ & 1
\end{bmatrix}\text{, and}
\\
&\mathcal C_{\not\square}\qquad \text{ represented by } \begin{bmatrix}
1 & u \\ & 1
\end{bmatrix}\text{, where $u\not\in \mathbb F_q^{\times 2}$.}
\end{split}
 \end{equation*}
Moreover, such an element $I - {}^\alpha X Y^{-1}X\in \mathcal C_\square$ (resp. $\mathcal C_{\not\square}$) if and only if $Y\in  \mathbb F_q^{\times 2}$ (resp. $Y\not\in  \mathbb F_q^{\times 2}$). The character values of $\rho_\pm$ at these classes are given by
$$\tr \rho_\pm \begin{bmatrix}
1 & u \\ & 1
\end{bmatrix}
=\frac{1}{2}\left({-1 \pm \omega(u)\mathfrak{n}(\psi) q^{1/2}}\right),\,\quad
u\in \mathbb F_q^\times,$$
where $\mathfrak{n}(\psi)$ is a 4th root of unity (the normalized Gauss sum depending on a chosen non-trivial additive character $\psi$ of $\mathbb F_q$) that we do not concern here.

Let $\zeta\in \mathbb F_q^\times \smallsetminus \mathbb F_q^{\times 2}$. When $\rho_y\neq \rho_z$, we derive from the character values that
$$\tr \rho_+ \begin{bmatrix}
1 & 1 \\ & 1
\end{bmatrix}
 - \tr \rho_+ \begin{bmatrix}
1 & \zeta \\ & 1
\end{bmatrix}
 = -\tr \rho_- \begin{bmatrix}
1 & 1 \\ & 1
\end{bmatrix}
+ \tr \rho_- \begin{bmatrix}
1 & \zeta \\ & 1
\end{bmatrix} = 2\mathfrak{n}(\psi)q^{1/2}.$$
The RHS of (\ref{traced_difference_Dy_Dz}) is then equal to 
$$ \left(\frac{q-1}{2}\right)(q^2-1)\left(2\mathfrak{n}(\psi)q^{1/2}\right),$$
 while the LHS is equal to 
$$( |b_y| - \nu |b_z|)\deg\rho_y = (1-\nu )q^{1/2}(q^2-1)\left(\frac{q-1}{2}\right),$$
so that we must have $\nu = -1$ (and also $A_y= \mathfrak n(\psi)$ in this case).

\addcontentsline{toc}{section}{References}

\bibliographystyle{alpha}
\bibliography{abc}

\end{document}